\documentclass[12pt,a4paper,reqno]{amsart}
\usepackage{amsmath}
\usepackage{amsfonts}
\usepackage{amssymb}
\numberwithin{equation}{section}

\theoremstyle{plain}

\newtheorem{theorem}[subsection]{Theorem}
\newtheorem{proposition}[subsection]{Proposition}
\newtheorem{lemma}[subsection]{Lemma}

\newtheorem{conjecture}[subsection]{Conjecture}

\newtheorem{remark}[subsection]{Remark}

\theoremstyle{definition}

\newtheorem{definition}[subsection]{Definition}

\renewcommand{\leq}{\leqslant}
\renewcommand{\geq}{\geqslant}

\newsavebox{\proofbox}
\savebox{\proofbox}{\begin{picture}(7,7)%
  \put(0,0){\framebox(7,7){}}\end{picture}}

\newcommand{\md}[1]{\ensuremath{(\mbox{mod}\, #1)}}
\newcommand{\mdsub}[1]{\ensuremath{(\mbox{\scriptsize mod}\, #1)}}

\def\B{{\mathcal B}}

\def\Z{{\mathbb Z}}
\def\E{{\mathbb E}}
\def\C{{\mathbb C}}
\def\R{{\mathbb R}}

\def\Q{{\mathbb Q}}

\def\eps{{\varepsilon}}

\def\emph#1{{\it #1}}
\def\textbf#1{{\bf #1}}

\begin{document}

\title[Obstructions to uniformity and primes]{Obstructions to uniformity, and arithmetic patterns in the primes}

\author{Terence Tao}
\address{Department of Mathematics\\University of California at Los Angeles\\ Los Angeles CA 90095}

\email{tao@math.ucla.edu}

\thanks{The author is
supported by a grant from the Packard Foundation.}

\subjclass{11N13, 11B25, 374A5}

\begin{abstract}  In this expository article, we describe the recent approach, motivated by ergodic theory, towards detecting arithmetic patterns in the primes, and in particular establishing in \cite{gt-primes} that the primes contain arbitrarily long arithmetic progressions.  One of the driving philosophies is to identify
precisely what the \emph{obstructions} could be that prevent the primes (or any other set) from behaving ``randomly'', and then either show
that the obstructions do not actually occur, or else convert the obstructions into usable structural information on the primes.
\end{abstract}

\maketitle

\section{Introduction}

An important class of problems in additive number theory, many of which are still far from being solved,
concerns the existence and distribution of affine-linear arithmetic patterns in the 
primes and almost primes.  Some well-known examples of these problems include:

\begin{itemize}
\item (Twin prime conjecture) Does there exist infinitely many numbers $n$ such that $n, n+2$ are both prime?
\item (Chen's theorem) \cite{Chen} There exists infinitely many numbers $n$ such that $n$ is prime, and $n+2$ is the product of at most two primes.
\item (Sophie Germain prime conjecture) Does there exist infinitely many numbers $n$ such that $n, 2n+1$ are both prime?
\item (Goldbach conjecture) For every sufficiently large even number $N$, does there exist an $n$ such that $n$ and $N-n$ are both prime?
\item (Vinogradov's theorem) \cite{vinogradov} For every sufficiently large odd number $N$, there exists $n,m$ such that $n$, $m$, and $N-n-m$ are all prime.
\item (Hardy-Littlewood prime tuples conjecture) \cite{hardy-littlewood} For any integers $a_1,\ldots,a_k$, which do not fill out all the residue classes of $\Z/p\Z$ for any prime $p$, there exists infinitely many $n$ such that $n+a_1,\ldots,n+a_k$ are all prime.
\item (van der Corput's theorem) \cite{van-der-corput} There exist infinitely many positive numbers $a,r$ such that $a, a+r, a+2r$ are all prime.
\item (Green-Tao theorem) \cite{green-tao} For any $k$, there exist infinitely many positive integers $a,r$ such that $a, a+r, \ldots, a+(k-1)r$ are all prime.
\end{itemize}

A unifying conjecture that encompasses all of these results is the \emph{generalized Hardy Littlewood prime tuples conjecture}, which
we now discuss.  As is customary in additive number theory, the most convenient way to count patterns in the primes is to introduce
the \emph{von Mangoldt function} $\Lambda: \Z \to \R^+$, defined by
setting $\Lambda(n) := \log p$ whenever $n = p^j$ is a power of a prime $p$ for some $j \geq 1$, and $\Lambda(n) = 0$ otherwise (in particular $\Lambda$ vanishes on zero and the negative integers). This function is mostly supported on the primes, and obeys a number of useful properties; for instance, one can encode the unique factorization of the integers
via the pleasant identity\footnote{All sums shall be over the positive integers $\Z^+$ unless otherwise indicated.}
\begin{equation}\label{logn}
 \log n = \sum_{d | n} \Lambda(d)
\end{equation}
for all $n \in \Z^+$.  Also, the prime number theorem can be phrased succinctly as
\begin{equation}\label{pnt}
\E( \Lambda(n) | 1 \leq n \leq N ) = 1 + o_{N \to \infty}(1)
\end{equation}
where we use $\E( f(n) | n \in A )$ to denote the average $\frac{1}{|A|} \sum_{n \in A} f(n)$, and $o_{N \to \infty}(1)$ denotes a quantity that goes to zero\footnote{Of course, one can make the decay rates much more quantitative, especially if one assumes strong hypotheses such as the Riemann hypothesis.  However, our discussion here will be not require any quantitative control of $o(1)$ type error terms.} as $N \to \infty$.
Thus $\Lambda$ is essentially normalized to have mean $1$.  More generally, for any modulus $q \geq 1$ and any integer $a$, we have
\begin{equation}\label{q-project}
\E( \Lambda(n) | 1 \leq n \leq N; n = a \md q ) = \Lambda_{\Z/q\Z}(a) + o_{N \to \infty;q}(1)
\end{equation}
for all sufficiently large $N$, where $o_{N \to \infty;q}(1)$ is a quantity which goes to zero as $N \to \infty$ for any fixed $q$, and the ``local von Mangoldt function'' $\Lambda_{\Z/q\Z}(a)$ is defined as the function which equals $\frac{q}{\phi(q)}$ when $a$ is 
coprime to $q$ and $0$ otherwise, with $\phi(q) = |(\Z/q\Z)^\times|$ being the Euler totient function; this result follows by combining the prime number theorem \eqref{pnt} with Dirichlet's theorem on the distribution of primes in arithmetic progressions.  One can also think of
\eqref{q-project} as an assertion that the $\Lambda_{\Z/q\Z}$ is essentially
the \emph{conditional expectation} of $\Lambda$ to the $\sigma$-algebra generated by the residue classes modulo $q$.

From the sieve of Eratosthenes, one is led to the heuristic\footnote{If $P$ is a statement, we use $1_P$ to denote the quantity $1$ if $P$ is true and $0$ if $P$ is false.  Similarly if $A$ is a set, we write $1_A(n)$ for $1_{n \in A}$.}
$$ \Lambda(n) \approx 1_{n > 0} \prod_{p < R} \Lambda_{\Z/p\Z}(n)$$
where $1 \ll R \ll n$ is an intermediate quantity between $1$ and $n$ that we shall be deliberately vague about specifying\footnote{The original sieve of Eratosthenes requires $R = \sqrt{n}$, but this is problematic for a number of reasons, for instance Mertens' theorem shows that a further correction term is required.  In practice we shall think of $R$ as being somewhat smaller, for instance a small power of $n$.}.  The Chinese remainder theorem
then suggests that the local factors $\Lambda_{\Z/p\Z}(n)$ in this product should behave ``independently''.  This leads to the following conjecture:

\begin{conjecture}[Generalized Hardy-Littlewood prime tuples conjecture]\label{hlconj}  Let $m, t$ be positive integers.  For each $1 \leq i \leq m$, let
$\psi_i: \Z^t \to \Z$ be an affine-linear form $\psi_i(x_1,\ldots,x_t) = \sum_{j=1}^t L_{ij} x_j + b_i$ for some integers $L_{ij}, b_i$, such that
the forms $\psi_i$ are all non-constant, and no two are rational multiples of each other.  Let $N$ be a large integer, and
assume that $b_i = O(N)$ for all $1 \leq i \leq m$.  Then we have
\begin{equation}\label{pll}
 \E( \prod_{i=1}^m \Lambda(\psi_i(x)) | x \in \{1,\ldots,N\}^t ) = \alpha_\infty(N) \prod_p \alpha_p + o_{N \to \infty; m,t,L}(1)
 \end{equation}
where $L := (L_{ij})_{1 \leq i \leq m, 1 \leq j \leq t}$, 
$\alpha_\infty(N)$ is the local density at infinity
$$ \alpha_\infty(N) := \E( \prod_{i=1}^m 1_{\psi_i(x) > 0} | x \in \{1,\ldots,N\}^t )$$
and $\alpha_p$ is the local density at each prime $p$
\begin{equation}\label{alphap-def}
 \alpha_p := \E( \prod_{i=1}^m \Lambda_{\Z/p\Z}(\psi_i(x)) | x \in (\Z/p\Z)^t ).
\end{equation}
\end{conjecture}

\begin{remark}
The density $\alpha_\infty(N)$ simply reflects the fact that the primes are positive; this factor is just $1$ if all the $L_{ij}$ and $b_i$ are positive. Note we allow the $b_i$ to depend on $N$, and the error term $o_{N \to \infty; m,t,L}(1)$ is presumed to be independent of the $b_i$;
this is necessary in order for this conjecture to encompass such conjectures as Goldbach's conjecture.
One can show that $\alpha_p = 1 + O_{m,t,L}(1/p^2)$ and hence the product $\prod_p \alpha_p$ (also known as the \emph{singular series})
is always convergent.  The conjecture is an assertion that the von Mangoldt function $\Lambda(n)$ behaves ``randomly'', subject to the structural constraints that it must resemble $1_{n > 0}$ ``locally at infinity'' (e.g. in the sense of \eqref{pnt}), and must 
resemble $\Lambda_{\Z/p\Z}$ locally at each prime $p$ (e.g. in the sense of \eqref{q-project}).  One can also extend the conjecture
to polynomial $\psi_i$; this is known as the Bateman-Horn conjecture \cite{bh}.
\end{remark}

This conjecture, if true, would imply all the conjectures and theorems stated earlier.  For instance, it predicts
\begin{equation}\label{lnn}
\E( \Lambda(n) \Lambda(n+2) | 1 \leq n \leq N ) = \prod_p \alpha_p + o_{N \to \infty}(1)
\end{equation}
where $\alpha_2 := 2$ and $\alpha_p := 1 - \frac{1}{(p-1)^2}$ for all odd primes $p$.  The \emph{twin prime constant} 
$$\Pi_2 := \prod_{p \hbox{ odd}} \alpha_p = 0.66016\ldots > 0$$
is positive, and \eqref{lnn} can then easily be seen to imply the twin prime conjecture.  Similarly for the other conjectures and theorems
stated earlier.

Of course, this conjecture is still hopelessly out of reach in the general case.  However, several partial results are known.  The bounds
\eqref{pnt}, \eqref{q-project} can already handle the $m=1$ case of this conjecture and more generally they can handle any ``non-degenerate''
case with $m \leq t$.  The Hardy-Littlewood circle method, which we discuss below, is roughly speaking able to handle any non-degenerate
case with $3 \leq m \leq t+1$ (thus encompassing Vinogradov's theorem and van der Corput's theorem), as well as a few additional cases\footnote{For instance, by a clever iteration of the circle method, it was established in \cite{balog1} that for any $k$ there exist infinitely many $k$-tuples of distinct primes $p_1,\ldots,p_k$, such that all the midpoints $(p_i + p_j)/2$ are also prime.}, but does not seem able to handle the general
case.  The conjecture is also known to be true
if one averages over a suitable subset of the parameters $L_{ij}$, $b_i$; see \cite{balog-hl}.  In the general case, the technique of 
upper bound sieves in sieve theory can usually yield an upper bound of $C_{m,t} \alpha_\infty(N) \prod_p \alpha_p + o_{N \to \infty; m,t,L}(1)$
for \eqref{pll} for some explicit $C_{m,t}$ (which usually has to be at least $2$, thanks to the notorious \emph{parity problem}); see also Section \ref{pnt-sec} below.
Closely related to this are the results of Goldston and Y{\i}ld{\i}r{\i}m, which show that asymptotic formulae such as \eqref{pll} can be 
recovered (but again with a loss of $C_{m,t}$ on the right-hand side) if one replaces $\Lambda$ with a slightly larger function $\nu$ 
which is localized to \emph{almost primes} (numbers with no small divisors) rather than primes themselves.
The ergodic theory-style transference arguments used in \cite{gt-primes}, \cite{green-tao} can conversely give \emph{lower bounds} of
$c_{m,t} \alpha_\infty(N) \prod_p \alpha_p + o_{N \to \infty; m,t,L}(1)$ for some small $0 < c_{m,t} < 1$, but only for linear forms
which are \emph{homogeneous} (no constant term $b_i$) and which are \emph{translation invariant}, in the sense that they take the form
$$ \psi_i(x_1,\ldots,x_t) = x_1 + \tilde \psi_i(x_2,\ldots,x_t).$$
In this special case, which covers the case of arithmetic progressions in the primes, there is also some hope of recovering the full
asymptotic \eqref{pll}; we discuss this below.

In this expository article we shall discuss these techniques, starting with the prime number theorem (but re-interpreted in the perspective
of Goldston-Y{\i}ld{\i}r{\i}m majorants), the 
classical circle method (but re-interpreted in a more ``ergodic'' perspective),  
and then turning to long arithmetic progressions in the primes; we also discuss some further recent progress in the case of progressions of length four.  In particular we hope to communicate some of the main philosophical ideas underlying 
the approach in \cite{gt-primes}, namely:

\begin{itemize}

\item Viewing the primes as a dense subset, not of the integers, but instead of a ``pseudorandom'' set of almost primes (or more precisely, a pseudorandom majorant $\nu$ for the von Mangoldt function $\Lambda$);
\item Attacking problems such as \eqref{pll} by locating the ``obstructions to uniformity'' which could potentially prevent \eqref{pll} from being true;
\item Using tools such as conditional expectation to handle these obstructions to uniformity, or tools such as the circle method to show that they do not occur at all.

\end{itemize}

This is by no means intended to be an exhaustive survey; see for instance \cite{kumchev} for a more in-depth discussion of many of these issues.
We will also not give detailed proofs for most of the assertions in this survey, referring the reader instead to the original papers.

\section{The prime number theorem and enveloping sieves}\label{pnt-sec}

We begin with the classical prime number theorem \eqref{pnt}.  The story of this theorem, and its connection to the zeroes of the Riemann zeta function
$\zeta(s) := \sum_n \frac{1}{n^s}$, is of course very well known, but we revisit it to make two points.  Firstly, as was observed by Chebyshev, one can obtain upper and lower bounds for \eqref{pnt} by elementary means (utilizing the pole of $\zeta$ at $s=1$, but requiring no further knowledge about zeroes or analytic continuation) that are only off by an absolute constant.  Secondly, by a refinement of this elementary method one can in fact get asymptotics with $o(1)$ error terms, but at the cost of smoothing out the von Mangoldt function $\Lambda$ and replacing it by a slightly larger variant, namely an \emph{enveloping sieve} $\nu$ for $\Lambda$.  In fact, it turns out even such results as those in \cite{gt-primes}, establishing arbitrarily long arithmetic progressions in the primes,
can in fact be proven without knowledge of the full prime number theorem (and thus without knowing any non-trivial zero-free region for $\zeta$, or for
any other $L$-function), instead using only\footnote{Of course, the larger the zero-free region is known for the zeta function, the better the bounds one will obtain on the number of progressions, but if one just wants to obtain the qualitative result that there are infinitely many progressions, no zero-free region beyond the trivial one used here is required.} these elementary techniques, albeit in conjunction with a deep and powerful theorem
of Szemer\'edi.

We begin with the argument of Chebyshev (rephrased here in modern language).  If $s$ is any complex number with $\Re(s) > 1$, we may multiply \eqref{logn} by $\frac{1}{n^s}$ and sum in $n$, and make the change of variables $n=dm$, to obtain
$$ \sum_n \frac{\log n}{n^s} = \sum_d \frac{\Lambda(d)}{d^s} \sum_m \frac{1}{m^s} = \sum_d \frac{\Lambda(d)}{d^s} \zeta(s).$$
The right-hand side is $-\zeta'(s)$, and hence we have the standard formula
\begin{equation}\label{lambda-form}
\sum_d \frac{\Lambda(d)}{d^s} = - \frac{\zeta'(s)}{\zeta(s)}.
\end{equation}
From summation by parts we obtain the bounds
\begin{equation}\label{zeta-control}
\zeta(s) = \frac{1}{s-1} + O(1); \quad \zeta'(s) = \frac{1}{(s-1)^2} + O(1)
\end{equation}
when $\Re(s) > 1$ and $s$ is close to 1.  In particular, we have a very small zero free region for $\zeta$ near $s=1$.
We conclude that
\begin{equation}\label{lambdadds}
 \sum_d \frac{\Lambda(d)}{d^s} = \frac{1}{s-1} + O(1)
 \end{equation}
whenever $\Re(s) > 1$ and $s$ is close to 1.  This, combined with the trivial observation that $\Lambda$ is non-negative, is already enough
to give the elementary bounds
\begin{equation}\label{crude}
 c - o_{N \to \infty}(1) \leq \E( \Lambda(n) | 1 \leq n \leq N ) \leq C + o_{N \to \infty}(1)
\end{equation}
for some absolute constants $0 < c < 1 < C$; for instance the upper bound follows by setting $s := 1 + \frac{1}{\log N}$ in \eqref{lambdadds}, while
the lower bound follows by setting $s := 1 + \frac{C'}{\log N}$ for some large $C'$ and using the upper bound already obtained to eliminate error terms.

The estimate \eqref{crude} is not an asymptotic, of course, since $c \neq C$.  However, we can recover good asymptotics by smoothing out
the von Mangoldt function $\Lambda$ slightly.  We introduce the \emph{M\"obius function} $\mu: \Z^+ \to \{-1,0,+1\}$, defined by $\mu(n) = (-1)^k$
when $n$ is the product of $k$ distinct primes for some $k \geq 0$, and $\mu(n) = 0$ otherwise.  The significance of this function lies in the inclusion-exclusion formula
\begin{equation}\label{mud-invert}
 1_{n=1} = 1_{n > 0} \sum_{d|n} \mu(d),
 \end{equation}
and hence from \eqref{logn}
\begin{equation}\label{lambdamu}
\begin{split}
\Lambda(n) &= 1_{n > 0} \sum_{m|n} \Lambda(m) 1_{n/m = 1} \\
&= 1_{n > 0} \sum_{dm|n} \Lambda(m) \mu(d) \\
&= 1_{n > 0} \sum_{d|n} \mu(d) \log \frac{n}{d} \\
&= 1_{n > 0} \log n \sum_{d|n} \mu(d) (1 - \frac{\log d}{\log n}).
\end{split}
\end{equation}
Inspired by this, let us define the truncated von Mangoldt functions $\Lambda_{\R,\varphi}: \Z \to \R$ by
\begin{equation}\label{lurch}
 \Lambda_{R,\varphi}(n) := \log R \sum_{d|n} \mu(d) \varphi( \frac{\log d}{\log R} )
\end{equation}
where $R > 1$ is a large parameter, and $\varphi: \R \to \R$ is a function supported on the interval $[-1,1]$.  For instance,
the von Mangoldt function itself corresponds to the case when $R = n$ and $\varphi(x) := \max(1-|x|,0)$.  The case
when $R < n$ and $\varphi(x) = \max(1-|x|,0)$ was studied by Goldston and Y{\i}ld{\i}rim; that case is also related to the Selberg upper bound sieve\footnote{The choice $\varphi(x) = \max(1-|x|,0)$ will give an optimized value of the relative density between $\Lambda$ and its enveloping sieve, although we will not need such optimization in our arguments.  Very recently, however, there has been work of Goldston, Motohashi, Pintz, and Y{\i}ld{\i}r{\i}m, which use precise optimization of higher-dimensional enveloping sieves in order to establish small gaps between primes, thus exploiting enveloping sieves in a rather different way than that discussed here.}, see \cite{green-tao} for further discussion.  These functions are more ``localized'', and hence easier to analyze, than the original von Mangoldt function, in the sense that they only involve divisors $d$ that are less than $R$\footnote{This can be viewed as a manifestation of the \emph{uncertainty principle}: localizing a function in the spectral or ``frequency'' sense (i.e. with respect to the divisors $d$)
must necessarily cause delocalization in physical space (i.e with respect to the variables $n$).}.  

The truncated von Mangoldt functions behave somewhat similarly to the von Mangoldt function, but are concentrated on the \emph{almost primes} rather than the primes themselves.  For instance, it is easy to see that $\Lambda_{R,\varphi}(n) = \varphi(0) \log R$ whenever $n$ is a prime larger than $R$, or more generally if $n$ is the product of primes larger than $R$.  One can also easily establish a fairly elementary 
``prime number theorem'' for these functions, provided that $R$ is not quite as large as $N$:

\begin{proposition}[Prime number theorem for $\Lambda_{R,\varphi}$]\label{pnt-lambda}  If $N^\eps \leq R \leq N^{1-\eps}$ for some $\eps > 0$, and $\varphi$ is smooth with $\varphi(0) = 1$ and $\varphi'(0) = 0$, then we have
\begin{equation}\label{lrchi}
\E( \Lambda_{R,\varphi}(n) | 1 \leq n \leq N ) = 1 + o_{N \to \infty; \eps,\varphi}(1).
\end{equation}
\end{proposition}

\begin{proof}  We can expand the left-hand side of \eqref{lrchi} as
$$ \log R \sum_{d \leq R} \mu(d) \varphi\left(\frac{\log d}{\log R}\right) \E( 1_{d|n} | 1 \leq n \leq N ).$$
From the elementary estimate 
$$ \E( 1_{d|n} | 1 \leq n \leq N ) = \frac{1}{d} + O( \frac{1}{N} )$$
we can thus write the left-hand side of \eqref{lrchi} as
$$ \log R \sum_{d \leq R} \frac{\mu(d)}{d} \varphi\left(\frac{\log d}{\log R}\right) + O_\varphi( \log R \sum_{d \leq R} \frac{1}{N} ).$$
Here the subscripting of $O()$ by $\varphi$ denotes that the implied constant is allowed to depend on $\varphi$.
Since $\varphi$ is supported on $[-1,1]$, we may remove the restriction $d \leq R$.
Since we are taking $R \leq N^{1-\eps}$, the error term here is $o_{N \to \infty; \eps,\varphi}(1)$.  Since we also take $R > N^\eps$, it
thus suffices to show that
\begin{equation}\label{lrchi-2}
 \log R \sum_{d} \frac{\mu(d)}{d} \varphi\left(\frac{\log d}{\log R}\right) =
1 + o_{R \to \infty; \varphi}(1).
\end{equation}
To proceed further we need to split $\varphi(\frac{\log d}{\log R})$ into expressions which are multiplicative in $d$.  This is easiest to establish by Fourier expansion\footnote{One could also use contour integration methods here instead of Fourier methods; the two approaches are essentially equivalent.}.  Since the function $e^x \varphi(x)$ is smooth and compactly supported, we have
\begin{equation}\label{ech}
 e^x \varphi(x) = \int_{-\infty}^\infty \psi(t) e^{-ixt}\ dt
 \end{equation}
for some rapidly decreasing function\footnote{In other words, $\psi(x) = O_{A,\psi}( (1 + |x|)^{-A} )$ for all $A > 0$ and $x \in \R$.}  
$\psi$.  We truncate this at $|t| = \log^{1/2} R$ (for instance) to obtain
$$ e^x \varphi(x) = \int_{|t| \leq \log^{1/2} R} \psi(t) e^{-ixt}\ dt + O_{A,\varphi}(\log^{-A} R)$$
for any $A > 0$.
In particular, we have
\begin{equation}\label{chir}
 \varphi(\frac{\log d}{\log R}) = \int_{|t| \leq \log^{1/2} R} \frac{\psi(t)\ dt}{d^{(1+it)/\log R}} + O_{A,\varphi}(d^{-1/\log R} \log^{-A} R)
 \end{equation}
and hence the left-hand side of \eqref{lrchi-2} can be written as
$$ \log R \int_{|t| \leq \log^{1/2} R} [\sum_d \frac{\mu(d)}{d^{1 + (1+it)/\log R}}] \psi(t)\ dt
+ O_{A,\varphi}( \log R \sum_d \frac{1}{d} d^{-1/\log R} \log^{-A} R ).$$
By taking $A=3$ (say), we see that the error term is $o_{R \to \infty;\varphi}(1)$ and so can be discarded.  As for the main term, we first
repeat the derivation of \eqref{lambda-form}, using \eqref{mud-invert} instead of \eqref{logn}, to conclude
$$ \sum_d \frac{\mu(d)}{d^s} = \frac{1}{\zeta(s)};$$
by \eqref{zeta-control} we thus have
$$ \sum_d \frac{\mu(d)}{d^s} = s-1 + O(|s-1|^2)$$
when $\Re(s) > 1$ and $s$ is sufficiently close to 1.  Setting $s = 1 + \frac{1+it}{\log R}$ for some $|t| \leq \log^{1/2} R$ 
we obtain (for $N$ and hence $R$
sufficiently large)
$$ \sum_d \frac{\mu(d)}{d^{1 + (1+it)/\log R}} = \frac{1+it}{\log R} + O( (1+|t|^2) \log^{-2} R ).$$
Inserting this bound into the previous computations, and using the rapid decay of $\psi$, we can thus write the left-hand side of \eqref{lrchi-2} as
$$ \int_{|t| \leq \log^{1/2} R} (1+it) \psi(t)\ dt
+ o_{R \to \infty;\varphi}(1).$$
Using the rapid decay of $\psi$ again, we can write this as
$$ \int_{-\infty}^\infty (1+it) \psi(t)\ dt
+ o_{R \to \infty;\varphi}(1)$$
which we rewrite in turn as
$$ (1 - \frac{d}{dx}) \int_{-\infty}^\infty e^{-ixt} \psi(t)\ dt|_{x=0}
+ o_{R \to \infty;\varphi}(1).$$
Applying \eqref{ech}, this becomes
$$ \varphi(0) - \varphi'(0) + o_{R \to \infty;\varphi}(1),$$
and the claim follows from the hypotheses on $\varphi$.
\end{proof}

One notable drawback of the truncated von Mangoldt functions $\Lambda_{R,\varphi}$ is that, unlike $\Lambda$, it is perfectly possible
for $\Lambda_{R,\varphi}(n)$ to be negative.  This however can be rectified by replacing $\Lambda_{R,\varphi}$ with the variant
\begin{equation}\label{nude}
\nu = \nu_{R,\varphi} := \frac{1}{\log R} \Lambda_{R,\varphi}^2.
\end{equation}
This function is still large on almost primes, indeed $\nu(n) = \Lambda_{R,\varphi}(n) = \varphi(0)^2 \log R$ whenever $n$ is a prime greater than $R$, or a product of primes greater than $R$.  In particular, if $\log R \sim \log N$ and $\varphi(0) \sim 1$ then we have the pointwise bound
\begin{equation}\label{lambda}
0 \leq \Lambda(n) \leq C \nu(n)
 \end{equation}
for all $1 \leq n \leq N$, where $C := \frac{1}{|\varphi(0)|^2} \frac{\log N}{\log R}$.
As observed\footnote{Strictly speaking, these authors only consider the case $\varphi(x) = \max(1-|x|,0)$, but the argument extends to general $\varphi$ without difficulty.} by Goldston and Y{\i}ld{\i}r{\i}m, we 
can also modify the above argument to obtain a prime number theorem for $\nu$, although at the cost of reducing the size of $R$:

\begin{proposition}[Prime number theorem for $\nu$]\label{pnt-nu}  
If $N^\eps \leq R \leq N^{1/2-\eps}$ for some $\eps > 0$, and $\varphi$ is smooth with $\int_0^1 |\varphi'(x)|^2\ dx = 1$, then we have
\begin{equation}\label{lrchi-nu}
\E( \nu(n) | 1 \leq n \leq N ) = 1 + o_{N \to \infty; \eps,\varphi}(1).
\end{equation}
\end{proposition}

\begin{proof} We repeat the proof of Proposition \ref{pnt-lambda}.  We can expand the left-hand side of \eqref{lrchi-nu} as
$$ \log R \sum_{d, d' \leq R} \mu(d) \mu(d') \varphi(\frac{\log d}{\log R}) \varphi(\frac{\log d'}{\log R}) 
\E( 1_{d,d'|n} | 1 \leq n \leq N ).$$
From the Chinese remainder theorem we have
$$ \E( 1_{d,d'|n} | 1 \leq n \leq N ) = \frac{1}{[d,d']} + O( \frac{1}{N} )$$
where $[d,d']$ is the least common multiple of $d$ and $d'$.  The hypothesis $R \leq N^{1/2-\eps}$ allows us to discard the error term as before, leaving us with the task of establishing
$$ \log R \sum_{d, d'} \frac{\mu(d) \mu(d')}{[d,d']} \varphi(\frac{\log d}{\log R}) \varphi(\frac{\log d'}{\log R}) 
= 1 + o_{R \to \infty; \varphi}(1).$$
From \eqref{chir} we have
$$
\varphi(\frac{\log d}{\log R}) \varphi(\frac{\log d'}{\log R}) = 
\int_{|t|, |t'| \leq \log^{1/2} R} \frac{\psi(t)\psi(t')\ dt dt'}{d^{(1+it)/\log R} (d')^{(1+it')/\log R}} + 
O_{A,\varphi}((dd')^{-1/\log R} \log^{-A} R).$$
Let us first dispose of the error term.  This contribution can be bounded by
$$ O_{A,\varphi}( \log R)^{1-A} \sum_{d,d'} \frac{1}{[d,d'] (dd')^{-1/\log R} }.$$
Using unique factorization $\Z^+ = \prod_p p^{\Z^+}$, and the multiplicative nature of the summand,
the sum can be expanded as an Euler product
$$ \sum_{d,d'} \frac{1}{[d,d'] (dd')^{-1/\log R} }
= \prod_p \sum_{d,d' \in p^{\Z^+}} \frac{1}{[d,d'] (dd')^{-1/\log R} }.$$
One can compute
$$ \sum_{d,d' \in p^{\Z^+}} \frac{1}{[d,d'] (dd')^{-1/\log R} } = 1 + O(1/p^{1+1/\log R}) \leq (1 - 1/p^{1+1/\log R})^{-O(1)}.$$
On the other hand, from \eqref{zeta-control} and the Euler product
$$ \zeta(s) = \prod_p \sum_{n \in p^{\Z^+}} \frac{1}{n^s} = \prod_p (1 - 1/p^s)^{-1}$$
we have
\begin{equation}\label{zeta-control-2}
\prod_p (1 - 1/p^s)^{-1} = \frac{1}{s-1} + O(1)
\end{equation}
for $\Re(s) > 1$ and $s$ close to $1$.  From this we see that the total contribution of the error term is
$O_{A,\varphi}( \log^{O(1)-A} R)$, which is acceptable since $A$ can be chosen to be large.

It remains to control the main term, which is
\begin{equation}\label{maint}
\log R \int_{|t|, |t'| \leq \log^{1/2} R} \left(\sum_{d,d'} \frac{\mu(d) \mu(d')}{[d,d'] d^{(1+it)/\log R} (d')^{(1+it')/\log R}}\right) \psi(t)\psi(t')\ dt dt'.
\end{equation}
The expression inside the parentheses can be expanded as an Euler product
$$ \prod_p \sum_{d,d'\in p^{\Z_+}} \frac{\mu(d) \mu(d')}{[d,d'] d^{(1+it)/\log R} (d')^{(1+it')/\log R}}$$
which one can compute as
$$ \prod_p ( 1 - \frac{1}{p^{1 + (1 + it)/\log R}} - \frac{1}{p^{1 + (1 + it')/\log R}} + \frac{1}{p^{1 + (2 + it + it')/\log R}} ).$$
After some Taylor expansion, we can write this as
\begin{equation}\label{brak}
 \prod_p \frac{( 1 - \frac{1}{p^{1 + (1 + it)/\log R}}) ( 1 - \frac{1}{p^{1 + (1 + it')/\log R}}) }{ 1 - \frac{1}{p^{1 + (2 + it+it')/\log R}} }
(1 + O( \frac{(1 + |t| + |t'|) \log p}{p^2 \log R} ) ).
\end{equation}
Since $\sum_p \frac{\log p}{p^2}$ is convergent, and $|t|, |t'| \leq \log^{1/2} R = o_{R \to \infty}(\log R)$, we have
$$ \prod_p (1 + O( \frac{(1 + |t| + |t'|) \log p}{p^2 \log R} ) ) = 1 + o_{R \to \infty}(1)$$
Applying \eqref{zeta-control-2}, we can thus write \eqref{brak} as
$$ (1 + o_{R \to \infty}(1) ) \log^{-1} R \frac{(1+it)(1+it')}{2 + it + it'}.$$
The contribution of the error term to \eqref{maint} is $o_{R \to \infty,\varphi}(1)$, thanks to the rapid decrease of $\psi$.  Hence we are left with
the expression
$$ \int_{|t|, |t'| \leq \log^{1/2} R} \frac{(1+it)(1+it')}{2 + it + it'} \psi(t)\psi(t')\ dt dt'$$
and by using the rapid decay of $\psi$ again, we see that we will be done as soon as we establish the identity
$$ \int_{-\infty}^\infty \int_{-\infty}^\infty \frac{(1+it)(1+it')}{2 + it + it'} \psi(t)\psi(t')\ dt dt' = 1.$$
Since
$$ \frac{1}{2+it + it'} = \int_0^\infty e^{-(2+it+it') x}\ dx = \int_0^\infty e^{-(1+it) x} e^{-(1+it') x}\ dx$$
the left-hand side can be written as
$$ \int_0^\infty ( \int_{-\infty}^\infty \psi(t) (1 + it) e^{-(1+it) x}\ dx )^2\ dt.$$
But by dividing \eqref{ech} by $e^x$ and then differentiating in $x$, we obtain
$$ \varphi'(x) = -\int_{-\infty}^\infty \psi(t) (1+it) e^{-ixt}\ dt$$
and the claim follows.
\end{proof}

It turns out that the above elementary argument is quite flexible, and can also give more sophisticated estimates for $\nu$, similar to \eqref{pll}.
Indeed we have

\begin{theorem}[Generalized Hardy-Littlewood prime tuples conjecture for $\nu$]\label{prime-corr}  Let $m, t$ be positive integers.  For each $1 \leq i \leq m$, let
$\psi_i: \Z^t \to \Z$ be an affine-linear form $\psi_i(x_1,\ldots,x_t) = \sum_{j=1}^t L_{ij} x_j + b_i$ for some integers $L_{ij}, b_i$, such that
the forms $\psi_i$ are all non-constant, and no two are rational multiples of each other.  Let $N$ be a large integer, and
assume that $b_i = O(N)$ for all $1 \leq i \leq m$.  
If $N^\eps \leq R \leq N^{1/2m-\eps}$ for some $\eps > 0$, and $\varphi$ is smooth with $\int_0^1 |\varphi'(x)|^2\ dx = 1$, then we have
\begin{equation}\label{pll-nu}
 \E( \prod_{i=1}^m \nu(\psi_i(x)) | x \in \{1,\ldots,N\}^t ) = \prod_p \alpha_p + o_{N \to \infty; m,t,L,\eps,\varphi}(1)
\end{equation}
where $\alpha_p$ was defined in \eqref{alphap-def}.
\end{theorem}

We will not prove this result here, but remark that the proof is a routine extension of that used to prove Proposition \ref{pnt-nu}, and
very similar results were proven in \cite{goldston-yildirim-old1}, \cite{goldston-yildirim-old2},
\cite{goldston-yildirim}, \cite{gt-primes}, \cite{tao-remark}.  
One can also obtain moment bounds for $\nu$ in terms of various multilinear 
integrals involving $\psi$; see \cite{goldston-yildirim-old1}, \cite{goldston-yildirim-old2},
\cite{goldston-yildirim} for some computations of this sort.  The density at infinity, $\alpha_\infty$, is missing, because $\nu$ extends to the negative integers as well as the positive ones.
Note that as the order $m$ of the correlation increases, the range of available $R$
decreases, so if we set $R$ equal to a fixed power of $N$, we only obtain correlations to finitely high order.

In the language of \cite{ramare}, \cite{ramare-ruzsa}, the function $C \nu$ appearing in \eqref{lambda} is an \emph{enveloping sieve} for
the von Mangoldt function $\Lambda$.  Results such as Theorem \ref{prime-corr} establish correlation estimates for this sieve, which in
turn automatically imply \emph{upper} bounds for expressions such as \eqref{pll} which are off by a constant $C_{m,t} > 1$; thus the enveloping sieve can be used as an upper bound sieve, though it has many other uses also, thanks in large part to correlation estimates
such as\footnote{By modifying the enveloping sieve slightly, one can also get some useful estimates on the Fourier coefficients of $\nu$, see \cite{green-tao}.  Of course, similar estimates are also known for the Fourier coefficients of $\Lambda$ itself, though the estimates for $\nu$ are simpler and do not require the theory of Siegel zeroes. In particular, the estimates are effective without requiring strong hypotheses such as GRH.} Theorem \ref{prime-corr}.  More advanced methods in sieve theory can of course be used to reduce this loss $C_{m,t}$, although the parity problem prevents one from removing this constant entirely by sieve-theoretic methods.

We have asserted earlier that $\nu$ is concentrated on the almost primes, which are coprime to all numbers less than $R$.  Let us provide some 
further evidence of this claim.  From \eqref{ech} we have
$$ \varphi( \frac{\log d}{\log R} ) = \int_{-\infty}^\infty \psi(t) d^{-(1+it)/\log R}\ dt$$
and hence by \eqref{lurch}
$$ \Lambda_{R,\varphi}(n) = \log R \int_{-\infty}^\infty \psi(t) \sum_{d|n} \mu(d) d^{-(1+it)/\log R}\ dt.$$
We can factorize the sum as an Euler product
$$ \sum_{d|n} \mu(d) d^{-(1+it)/\log R}\ dt = \sum_{p|n} (1 - p^{-(1+it)/\log R})$$
and conclude
$$ \Lambda_{R,\varphi}(n) = \log R \int_{-\infty}^\infty \psi(t) \prod_{p|n} (1 - p^{-(1+it)/\log R})\ dt$$
and similarly by \eqref{nude}
$$ \nu(n) = \log R \int_{-\infty}^\infty \int_{-\infty}^\infty \psi(t) \psi(t') 
\prod_{p|n} (1 - p^{-(1+it)/\log R}) (1 - p^{-(1+it')/\log R})\ dt dt'.$$
Since $\psi(t)$ is rapidly decreasing, the integral effectively localizes $t$ to be close to 1.  The factor $(1 - p^{-(1+it)/\log R})$
is then close to $0$ when $p \ll R$ and oscillates around $1$ when $p \gg R$.  Thus we expect $\Lambda_{R,\varphi}(n)$ and $\nu(n)$
to be small when $n$ has one or more prime factors $\ll R$, and these quantities should 
be close to $\log R$ when $n$ is a product of primes $\gg R$, though in some exceptional cases
(when the phases of $p^{-(1+it)/\log R}$ align in an unfavourable way) one may expect $\Lambda_{R,\varphi}(n)$ to be somewhat larger than 
this\footnote{On the other hand, \eqref{lurch} shows that $\Lambda_{R,\varphi}(n)$ can be crudely bounded by $O_\varphi(\tau(n) \log R)$, where
$\tau(n) = \sum_{d|n} 1$ is the divisor function.  As is well known, the divisor function has size $O(\log n)$ on the average, though it can get
significantly larger than this for very smooth $n$.  However, it is always $O_\eps(n^\eps)$ for any $\eps > 0$, and hence $\Lambda_{R,\varphi}$ and $\nu$ also have this type of bound.}.  Thus we have the rough heuristics
\begin{equation}\label{apr}
\Lambda(n) \approx (\log N) 1_P; \quad \nu(n) \approx (\log R) 1_{AP}
\end{equation}
for $n \sim N$, where $P$ denotes the primes up to $N$ and $AP$ denotes the almost primes at level $R$ up to $N$ (i.e. the products of primes larger than $R$).  Observe that $\Lambda(n)$ and $\nu(n)$ both have average $1 + o_{N \to \infty}(1)$, which thus suggests that $P$ has density about
$\frac{\log R}{\log N}$ inside $AP$; one can obtain more precise estimates here using Buchstab's formula.  On the other hand, Theorem
\ref{prime-corr} combined with the heuristic \eqref{apr} suggests that the set $AP$ is very nicely distributed if $R = N^\eps$ for some suitably small $\eps$.  Thus, in summary, the primes $P$ form a set of positive density ($\approx \eps$) inside the almost primes at level $R = N^\eps$, and
the latter set has a well-controlled distribution.  This turns out to be a very useful perspective for a number of problems, as it bypasses
the difficulty that the primes have only a density of $\frac{1}{\log N} = o_{N \to \infty}(1)$ with respect to the integers $\{1,\ldots,N\}$.
Thus the almost primes $AP$ (or more precisely the enveloping sieve $\nu$) forms a better majorant for the primes $P$ (or more precisely
the von Mangoldt function $\Lambda$) than the integers (or the constant $1$).

\section{The $W$-trick}

As we have seen, any correlation estimate involving $\Lambda$ or $\nu$ will involve a number of local densities $\alpha_p$; these densities ultimately arise from the fact that the projections $\Lambda_{\Z/p\Z}$ of $\Lambda$ to the residue classes modulo $p$ are not constant.  Note that due to the rapid convergence of the product $\prod_p \alpha_p$, it is only the small divisors $p$ for which this non-uniformity is significant.
However, if one does not care much about the exact order of decay in the $o(1)$ errors, then there is a cheap trick, which we call the ``$W$-trick'', available to essentially eliminate the role of these local factors, so that one only has to deal with functions which are uniform with respect to small divisors.

This trick works as follows.  We introduce a new parameter $1 \ll w \ll N$; this will eventually be set to a very slowly growing function of $N$,
such as $\log \log N$, although for the purposes of getting qualitative $o(1)$ bounds it is not particularly important what $w$ is.  We let
$W := \prod_{p < w} p$ be the product of all the primes less than $w$.  The prime numbers larger than $w$ will then be distributed in the residue
classes $\{ Wn + b: n \in \Z\}$, where $b$ is one of the $\phi(W)$ numbers in $\{1,\ldots,W\}$ which are coprime to $W$.  For each of these numbers
$b$, we introduce the renormalized von Mangoldt function
$$ \Lambda_{b \mdsub W}(n) := \frac{W}{\phi(W)} \Lambda(Wn + b)$$
and similarly the renormalized truncated von Mangoldt functions
$$ \Lambda_{R,\varphi,b \mdsub W}(n) := \frac{W}{\phi(W)} \Lambda_{R,\varphi}(Wn + b)$$
and the renormalized enveloping sieve
$$ \nu_{b \mdsub W}(n) := \frac{W}{\phi(W)} \nu(Wn+b).$$
Then the functions $\Lambda_{b \mdsub W}(n)$ behave like $\Lambda$ except that the projections modulo $q$ are now extremely close to $1$
for small $q$.  Indeed from \eqref{q-project} and the Chinese remainder theorem, one easily verifies
$$
\E( \Lambda_{b \mdsub W}(n) | 1 \leq n \leq N; n = a \md q ) = 1 + o_{N \to \infty;w}(1)$$
for $1 \leq q \leq w$.  The analogue of Conjecture \ref{hlconj} is then the assertion that
\begin{equation}\label{pll-bi}
 \E( \prod_{i=1}^m \Lambda_{b_i}(\psi_i(x)) | x \in \{1,\ldots,N\}^t ) = \alpha_\infty(N) \prod_{p > w} \alpha_p + o_{N \to \infty; m,t,L,w}(1)
 \end{equation}
whenever $b_1,\ldots,b_m \in \{1,\ldots,W\}$ are coprime to $W$; thus the local factors corresponding to primes less than or equal to $w$
in \eqref{pll} are eliminated, at the cost of letting the $o(1)$ term depend on $w$.  
Actually it is not hard to see that \eqref{pll} is in fact equivalent to \eqref{pll-bi}.  In many cases, the remaining local
factor $\prod_{p > w} \alpha_p$ is in fact $1 + o_{w \to \infty; m,t,L}(1)$; for instance, this is the case if no two of the
linear parts $(L_{ij})_{1 \leq j \leq t}$ of the affine forms $\psi_1,\ldots,\psi_m$ are not rational multiples of each other.  However, there are some
important cases where the remaining local factors are significant.  For instance, for $1 \leq a \leq N$, the prime tuples conjecture
predicts
$$ \E( \Lambda_b(x) \Lambda_b(x+a) | 1 \leq x \leq N ) = \prod_{p > w: p | a} (1 + \frac{1}{p}) ( 1 + o_{w \to \infty}(1) + o_{N \to \infty;w}(1) ).$$
The expression $\tau(a) := \prod_{p > w: p | a} (1 + \frac{1}{p})$ is small for most $a$, for instance one can establish the moment estimates
\begin{equation}\label{tauq}
 \E( \tau^q(a) | 1 \leq a \leq N ) = O_q(1)
 \end{equation}
for all $1 \leq q < \infty$ (indeed one can refine the right-hand side to $1 + o_{w \to \infty;q}(1)$).  However it is not bounded, as can be seen by
taking $n$ to be the product of a large number of primes, each of which is slightly larger than $w$.  Nevertheless it is a good heuristic to
view $\prod_{p > w} \alpha_p$ as being close to $1$ for ``most'' choices of forms $\psi_i$.

Similar considerations apply to the enveloping sieve $\nu$.  For instance, one can establish that
\begin{equation}\label{nub}
\E( \prod_{i=1}^m \nu_{b_i}(\psi_i(x)) | x \in \{1,\ldots,N\}^t ) = 1 + o_{N \to \infty; m,t,L,\eps,\varphi,w}(1)
\end{equation}
whenever no two linear parts of the affine forms $\psi_1,\ldots,\psi_m$; this is essentially\footnote{The conditions verified in \cite{gt-primes} actually refer to a version of $\nu_b$ adapted to $\Z/N\Z$ rather than $\{1,\ldots,N\}$, but the distinction between the two is rather minor.} the \emph{linear forms condition}
verified in \cite[Proposition 9.8]{gt-primes}.  Similarly, one
can show\footnote{The diagonal cases
$a_i = a_j$ can be treated using the crude bound $\nu(n) = O_\eps(N^\eps)$ for any $\eps > 0$ and $n = O(N)$.}
\begin{equation}\label{corrtau}
 \E( \prod_{i=1}^m \nu_{b_i}(x + a_i) | x \in \{1,\ldots,N\}^t ) \leq \sum_{1 \leq i < j \leq m} \tilde \tau(a_i - a_j)
 \end{equation}
where $\tilde \tau: \Z \to \R^+$ is a slight variant of $\tau$ which is even and obeys the moment conditions \eqref{tauq}; this is essentially the \emph{correlation condition} verified in \cite[Proposition 9.10]{gt-primes}.  Morally, one should think of the right-hand side of \eqref{corrtau}
as being bounded, with only a few exceptions such as when $a_i - a_j$ is zero or very smooth (contains a large number of prime factors larger than $w$).

The linear forms condition \eqref{nub} is an assertion that the $\nu_b$ are distributed \emph{pseudorandomly} throughout $\{1,\ldots,N\}$; more informally,
the almost primes $AP$ when restricted to a coset $\{ Wn+b: n \in \Z \}$ with $b$ coprime to $W$, behave pseudorandomly inside each such coset.
This is consistent with the heuristics used to support the Hardy-Littlewood prime tuples conjecture, such as Cramer's probabilistic model for the primes.  In this context, a useful probabilistic model for $\nu_b(n)$ would be a function which equalled $\frac{W}{\phi(W)} \log R$ with probability
$\frac{W}{\phi(W) \log R}$ independently for each $n$, and equalled 0 otherwise.  The prime tuples conjecture then asserts that the $\Lambda_b$ also
behave in a similarly pseudorandom manner (but with $\log R$ essentially replaced by $\log N$).

The linear forms condition \eqref{nub} shows that the correlations of $\nu_b$ are very close to the correlations of the constant function $1$, thus
$\nu_b$ is close to $1$ in a ``weak'' sense.  One of the philosophies underlying the work in \cite{gt-primes} is a \emph{transference principle} which asserts, informally, that many results which are true for functions bounded by constant function $1$, are likely to extend to functions
bounded by pseudorandom functions such as $\nu_b$, or variants such as $\nu_b + 1$.

Any counting problem concerning the von Mangoldt function $\Lambda$ can of course be subdivided into a counting problem involving the $\Lambda_b$.
For instance, suppose one wanted to establish a bound such as
$$ \E( \Lambda(a) \ldots \Lambda(a+(k-1)r) | 1 \leq a,r \leq N ) \geq c_k - o_{N \to \infty;k}(1)$$
for all $k \geq 1$ and $N \geq 1$, and $c_k > 0$; this bound is in fact obtained in \cite{gt-primes}, and implies
that the primes contain arbitrarily long arithmetic progressions.  In order to achieve this bound, it suffices to show that for all $w$ there exist $b \in \{1,\ldots,W\}$ coprime to $W$ such that
\begin{equation}\label{lbr}
 \E( \Lambda_{b}(a) \ldots \Lambda_{b}(a+(k-1)r) | 1 \leq a,r \leq N ) \geq c'_k - o_{w \to \infty; k}(1) - o_{N \to \infty;k,w}(1)
 \end{equation}
for some other $c'_k > 0$.  Indeed,if such a bound were true, it would imply that
$$ \E( \Lambda_{b}(a) \ldots \Lambda_{b}(a+(k-1)r) | 1 \leq a,r \leq N ) \geq c'_k/2$$
(say) whenever $w$ was sufficiently large depending on $k$, and $N$ was sufficiently large depending on $w$ and $k$.  But
since $\Lambda_{b}$ is a renormalized component of $\Lambda$ using the affine-linear transformation $n \mapsto Wn+b$ (which preserves arithmetic
progressions), we then observe that
$$ \E( \Lambda(a) \ldots \Lambda(a+(k-1)r) | 1 \leq a,r \leq N ) \geq c_{k,w}$$
for some $c_{k,w} > 0$.  Fixing $w$ to be a suitably large constant depending only on $k$, we obtain the claim.

This reduction from $\Lambda$ to $\Lambda_b$ is used in \cite{gt-primes}.  Indeed, \eqref{lbr} is established for all $1 \leq b < W$ which are 
coprime to $W$.  In the proof, the only facts needed are the bounds $0 \leq \Lambda_b \leq C \nu_b$ (which is inherited from \eqref{lambda}) and
$\E( \Lambda_b(n) | 1 \leq n \leq N ) > c - o_{N \to \infty;w}(1)$ (which comes from \eqref{q-project}).  In fact, since we only need to establish
\eqref{lbr} for a \emph{single} $b$, it is possible to avoid using Dirichlet's theorem altogether, and simply use the pigeonhole principle to locate
a $b$ for which $\Lambda_b(n)$ has large mean.  This observation has the interesting application that it allows one to extend the result in \cite{gt-primes} to obtain arbitrarily long progressions, not just in the primes, but in fact in any subset of the primes (or almost primes)
of positive relative density.

In summary, the $W$-trick allows one to easily eliminate the influence of small divisors, resulting in functions $\Lambda_b$, $\nu_b$ which are much more uniformly distributed than their non-renormalized counterparts $\Lambda$, $\nu$.  Of course, the price one pays for doing so is that the $o(1)$ error terms, as well as the $c_k$ bounds employed above, deteriorate rather substantially; however if one is only interested in qualitative results then this trick is essentially cost-free.

\section{Fourier obstructions to uniformity}\label{fourier-sec}

We now discuss the problem of counting the progressions of length three in the primes.  This can of course be done by the circle method, and this is essentially what we do here, but we shall adopt the philosophy of counting progressions by first establishing what the obstructions are to uniformity,
and then dealing with these obstructions in some manner. The $W$-trick is already one way to eliminate one obstruction to uniformity, namely
irregular distribution when localized to small primes, which in the language of the circle method allows one to ignore the contribution of
the major arcs (except the major arc near $1$).  We will see other ways to deal with obstructions to uniformity later in this article.

The standard way to count progressions of length three in the primes is to try to obtain asymptotics, or at least bounds, for the average
\begin{equation}\label{lar}
\E( \Lambda(a) \Lambda(a+r) \Lambda(a+2r) | 1 \leq a,r \leq N ).
\end{equation}
Indeed Conjecture \ref{hlconj} already predicts an explicit asymptotic for this quantity, and Theorem \ref{prime-corr} gives an upper bound which
is only off by an absolute constant.  One would then use the Fourier transform right away, to convert this expression to an integral involving
an exponential sum such as $\E( \Lambda(n) e(-n\alpha) | 1 \leq n \leq N)$, where $\alpha$ is a real number and $e(x) := e^{2\pi i x}$.
This sum would then be estimated in two different ways, one when $\alpha$ is major arc (close to a rational with small denominator) and one
when $\alpha$ is minor arc.  The minor arc computation is reasonably elementary (ultimately relying on variants of the identity \eqref{lambdamu},
the Cauchy-Schwarz inequality, and some bilinear cancellation in the expression $e(-nm\alpha)$) but the major arc computation is somewhat deeper,
relying among other things on the Siegel-Walfisz theorem.

It turns out that one can proceed in a more elementary fashion if one is not seeking an asymptotic, but only a non-zero lower bound
on the quantity \eqref{lar} (which will certainly be enough to imply the qualitative result that there are infinitely many progressions of length
three in the primes).  Instead of needing to control the exponential sums of $\Lambda$, one only needs to control the exponential sums of a majorant
$\nu$ or $\nu_b$, which is much simpler.  However, one does need one additional ingredient, namely \emph{Roth's theorem} \cite{roth}.  Roth's original formulation of this theorem asserts that any subset of the integers with positive upper density, necessarily contains infinitely many progressions of length three.  Varnavides \cite{varnavides} showed that this qualitative version is in fact equivalent to the following more quantitative statement:

\begin{theorem}[Quantitative Roth theorem]\label{qrt}\cite{roth},\cite{varnavides}  Let $f: \Z/N\Z \to \R$ be a function such that $0 \leq f(n) \leq 1$ for all $n \in \Z/N\Z$, and such that $\E( f(n) | n \in \Z/N\Z ) \geq \delta$ for some $0 < \delta < 1$.  Then we have
$$ \E( f(a) f(a+r) f(a+2r) | a,r \in \Z/N\Z ) \geq c(\delta)$$
for some $c(\delta) > 0$.
\end{theorem}

The best value of $c(\delta)$ currently known is $c(\delta) \gg \delta^{C/\delta^2}$ for some absolute constant $C$, see
\cite{bourgain-triples}.  However for the qualitative arguments we give below, we do not need to know the exact value of $c(\delta)$.  We also do not
need to know the proof of Theorem \ref{qrt}; we may treat it as a ``black box''.  We do remark however that the known proofs of this theorem,
involving either Fourier analysis, ergodic theory, or graph theory, are extremely instructive and are very consistent with the philosophy outlined
here of detecting obstructions to uniformity and then somehow dealing with each of the obstructions which occur.  For us, the power of Roth's theorem 
lies in the fact that very little structural information is demanded of $f$ (in particular, no arithmetic structure or Fourier-analytic
structure is required), besides the important constraint that $f$ is bounded\footnote{Indeed, our entire philosophy here is in some sense the polar opposite of the more conventional approach, in which one builds up as much information about the primes (or any other number-theoretic object) as possible, for instance using deep estimates on Dirichlet $L$-functions, and then uses all this information to then attack quantities such as \eqref{lar}.  In contrast, we adopt a minimalist approach (in the spirit of sieve theory) in which we treat the primes as nothing more than a generic subset of the almost primes with positive relative density, ignoring all the rich arithmetic structure.  That this approach works at all, is entirely due to the existence of such
theorems as Roth's theorem, which apply to all sets of positive density (or bounded functions with large mean).  However, as we shall see later it is possible to blend the two approaches and use deeper facts about the primes to obtain sharper results.}.  

At present, Roth's theorem does not directly allow us to obtain any non-trivial lower bound on \eqref{lar} for two reasons.  The first (rather trivial) reason is that we have stated Roth's theorem in $\Z/N\Z$ rather than on $\{1,\ldots,N\}$, but there are some easy truncation arguments (which we omit) to pass back and forth between these two settings, possibly after modifying $N$ by a factor of $2$ or so.  The more serious difficulty is that
$\Lambda$ is not bounded, and if we do normalize $\Lambda$ to be bounded (e.g. by dividing by $\log N$) then $\delta$ becomes too small for Roth's theorem to be of any use.  However, as it turns out it is relatively easy to \emph{decompose} $\Lambda$ into a bounded function (for which Roth's theorem is applicable) and a ``uniform'' error (which has a negligible impact on \eqref{lar}).

Before we do this, we need to understand exactly what type of functions will give a negligible impact to expressions such as \eqref{lar}.  To phrase things a little more concretely, let us work in the cyclic group $\Z/N\Z$ instead of the progression $\{1,\ldots,N\}$, taking $N$ to be odd,
and consider an expression such as
\begin{equation}\label{fgh}
 \E( f(a) g(a+r) h(a+2r) | a,r \in \Z/N\Z )
\end{equation}
for some functions $f, g, h: \Z/N\Z \to \C$.  To begin the discussion let us take $f,g,h$ to be bounded in magnitude by $1$, although for applications
to the primes we will eventually need to discard this hypothesis.

Since $f,g,h$ are bounded by $1$, it is clear that \eqref{fgh} is also bounded in magnitude by $1$.  However, in many cases, \eqref{fgh} will
be much smaller than $1$.  For instance, if one of $f,g,h$ is small in some averaged sense, say if the $L^1$ norm $\E( |f(n)| | n \in \Z/N\Z )$ is small, then \eqref{fgh} will be small also.  Also, if one of $f, g, h$ fluctuates randomly, for instance if $f(n) = \pm 1$ for each $n$, with each
$f(n)$ attaining $+1$ or $-1$ independently with equal probability, then it is easy to see that \eqref{fgh} will be quite small with high probability.
Let us informally call a function \emph{linearly uniform}\footnote{The notation here is due to Gowers \cite{gowers-4}.  The term ``uniform'' arises because linearly uniform functions behave like a signed probabilistic point process with the uniform distribution; another possible terminology is ``linearly unbiased''.  Somewhat confusingly, this usage of the word ``uniform'' is completely different from, and in fact in opposition to, the notion of ``uniformly bounded''; indeed, we will later need to rely crucially on the fact that linearly uniform functions can be very far from being uniformly bounded.} if the 
expression \eqref{fgh} is necessarily small as soon as at least one of $f,g,h$
is set equal to this function.  Thus for instance functions with small $L^1$ norm, or randomly fluctuating functions, will be linearly uniform.
Since \eqref{fgh} is linear in $f$, $g$, and $h$ separately, we thus see that we can modify $f$, $g$, or $h$ by a linearly uniform function without
significantly affecting \eqref{fgh}, and so linearly uniform functions are ``negligible'' for the purposes of counting progressions of length
three.  On the other hand, from the identity
$$ \E( e(\alpha a) e(-2\alpha(a+r)) e(\alpha(a+2r)) | a,r \in \Z/N\Z) = 1$$
for any $\alpha \in \frac{1}{N} \Z$, we see that the function $n \mapsto e(\alpha n)$ is not linearly uniform.  More generally, since
$$ \E( f(a) e(-2\alpha(a+r)) e(\alpha(a+2r)) | a,r \in \Z/N\Z) = \E( f(n) e(-\alpha n) | n \in \Z/N\Z )$$
we see that any function $f$ which has a large correlation (inner product) with a linear phase function $e(\alpha n)$, will not be linearly uniform.
Thus linear phase functions are \emph{obstructions} to linear
uniformity; this may help explain the ``linear'' in the terminology ``linear uniformity''.

The effectiveness of the circle method, at least for the task of counting progressions of length three, ultimately lies in the fact that
linear phase functions are the \emph{only} obstructions to linear uniformity, at least when everything is bounded; thus if a bounded function has small correlation with every linear phase function,
then it is linearly uniform.  More precisely:

\begin{lemma}\label{fgh-lin}  Let $f,g,h: \Z/N\Z \to \C$ be functions bounded by $1$, and suppose that
$$ |\E( f(n) e(-\xi n / N) | n \in \Z/N\Z)| \leq \eps$$
for some $\eps > 0$ and all $\xi \in \Z/N\Z$.  Then we have
$$ | \E( f(a) g(a+r) h(a+2r) | a,r \in \Z/N\Z ) | \leq \eps.$$
\end{lemma}

Not co-incidentally, Lemma \ref{fgh-lin} is also the first step used in the Fourier-analytic proof of Roth's theorem; however, we will not discuss
this connection here.

\begin{proof}  Writing $\hat f(\xi) := \E( f(n) e(-\xi n/N) | n \in \Z/N\Z)$ for all $\alpha \in \Z/N\Z$, and similarly for
$\hat g$ and $\hat h$, we have the Fourier inversion formulae
\begin{align*}
f(a) &= \sum_{\xi \in \Z/N\Z} \hat f(\xi) e(\xi a/N); \\
g(a+r) &= \sum_{\lambda \in \Z/N\Z} \hat g(\lambda) e(\lambda(a+r)/N); \\
h(a+2r) &= \sum_{\eta \in \Z/N\Z} \hat h(\eta) e(\eta(a+2r)/N).
\end{align*}
Substituting these formulae and simplifying, we eventually obtain the identity
\begin{equation}\label{fgh-fourier}
\E( f(a) g(a+r) h(a+2r) | a,r \in \Z/N\Z ) = \sum_{\xi \in \Z/N\Z} \hat f(\xi) \hat g(-2\xi) \hat h(\xi).
\end{equation}
On the other hand, from Plancherel's identity and the boundedness of $g$ and $h$ we have
$$ \sum_{\xi \in \Z/N\Z} |\hat g(-2\xi)|^2 \leq 1; \quad \sum_{\xi \in \Z/N\Z} |\hat h(\xi)|^2 \leq 1$$
while from the hypothesis on $f$ we have $|\hat f(\xi)| \leq \eps$ for all $\xi$.  The claim then follows from H\"older's inequality.
\end{proof}

Now we return to the task of estimating \eqref{lar}.  Applying the $W$-trick to make $\Lambda$ more uniformly distributed, 
it will suffice to obtain an estimate of the form
$$ \E( \Lambda_b(a) \Lambda_b(a+r) \Lambda_b(a+2r) | 1 \leq a,r \leq N ) \geq c - o_{W \to \infty}(1) - o_{N \to \infty;W}(1)$$
for some absolute constant $c > 0$.  Let us cheat a little bit by identifying $\{1,\ldots,N\}$ with $\Z/N\Z$ (ignoring issues of truncation
and wraparound, which are actually not difficult to deal with), so that we are now faced with establishing a lower bound for
\begin{equation}\label{lab3}
\E( \Lambda_b(a) \Lambda_b(a+r) \Lambda_b(a+2r) | a,r \in \Z/N\Z ).
\end{equation}
We would like to use Lemma \ref{fgh-lin} to strip away the linearly uniform components of $\Lambda_b$.  However, we are faced with the difficulty
that $\Lambda_b$ is not uniformly bounded.  Fortunately, we can use the fact that $\Lambda_b$ is majorized by an enveloping sieve $\nu_b$.  Actually
we will not quite use the enveloping sieve $\nu_b$ constructed in the previous section, but use a slight variant $\tilde \nu_b$ which is closely
related to the Selberg sieve.  The enveloping sieve $\nu_b$ can be written down explicitly, but it is a little messy;
see \cite{green-tao} for a definition, together with a full analysis and comparison of these two enveloping sieves.  For this expository paper, suffice it to say that we still have the basic majorization
\begin{equation}\label{tnub}
0 \leq \Lambda_b \leq C \tilde \nu_b
\end{equation}
and that the Fourier coefficients of the Selberg enveloping sieve $\tilde \nu_b$ can be computed very explicitly; for instance one can show that
\begin{equation}\label{tnub-fourier}
\widehat{\tilde \nu_b}(\xi) = o_{W \to \infty}(1) + o_{N \to \infty;W}(1)
\end{equation}
for all $\xi \in \Z/N\Z \backslash \{0\}$.
Using this and other bounds, together with
orthogonality arguments such as those used in the large sieve (or of Tomas-Stein restriction theory), it is possible to obtain a weighted form
of the Plancherel theorem, namely that
\begin{equation}\label{lpp}
\| \hat f \|_{l^p(\Z/N\Z)} \ll_p 1 
\end{equation}
whenever $p > 2$ and $f: \Z/N\Z \to \C$ is bounded pointwise by $\tilde \nu_b + 1$; see \cite{green-tao} (and also \cite{green}).
The key point in these estimates is that no factor of $\log N$ appears on the right-hand side, despite the fact that all the $L^q$ moments of
$\Lambda$ and $\nu$ (except the $L^1$ moment) contains such a logarithmic factor.
Using this estimate we can obtain a weighted variant of Lemma \ref{fgh-lin}:

\begin{lemma}\label{fgh-lin-weighted}\cite{green-tao}  Let $f,g,h: \Z/N\Z \to \C$ be functions bounded in magnitude
by $\tilde \nu_b + 1$, and suppose that
$$ |\E( f(n) e(-\xi n/N) | n \in \Z/N\Z)| \leq \eps$$
for some $\eps > 0$ and all $\xi \in \Z/N\Z$.  Then we have
$$ | \E( f(a) g(a+r) h(a+2r) | a,r \in \Z/N\Z ) | \ll \eps^{1/2}.$$
\end{lemma}

\begin{proof}  From \eqref{fgh-fourier} and H\"older's inequality we have
$$ | \E( f(a) g(a+r) h(a+2r) | a,r \in \Z/N\Z ) | \leq \| \hat f \|_{l^\infty(\Z/N\Z)}^{1/2} \| \hat f \|_{l^{5/2}(\Z/N\Z)}^{1/2} \| \hat g \|_{l^{5/2}(\Z/N\Z)} \| \hat h \|_{l^{5/2}(\Z/N\Z)}$$
(for instance).  From hypothesis we have $\| \hat f \|_{l^\infty(\Z/N\Z)} \leq \eps$.  The claim now follows from \eqref{lpp}.
\end{proof}

Thus, even when considering functions that are merely bounded by $\tilde \nu_b + 1$ instead of bounded by $1$, it is still the case that linear phase
functions are the only obstruction to orthogonality.  One can view this as a weak version of Plancherel's theorem, transferred to the enveloping sieve $\tilde \nu_b + 1$.

At this point one could try to show that $\Lambda_b$, or more precisely the normalized function $\Lambda_b - 1$, has small correlation with all linear phase functions,
$$ \E( (\Lambda_b(n)-1) e(-\xi n / N) | n \in\Z/N\Z) = o_{W \to \infty}(1) + o_{N \to \infty;W}(1).$$
This, together with Lemma \ref{fgh-lin-weighted}, would imply that $\Lambda_b$ can be replaced with $1$ with negligible error in \eqref{lab3}
and we would conclude that
$$ \E( \Lambda_b(a) \Lambda_b(a+r) \Lambda_b(a+2r) | a,r \in \Z/N\Z ) = 1 + o_{W \to \infty}(1) + o_{N \to \infty;W}(1),$$
which would of course be consistent with the Hardy-Littlewood prime tuples conjecture.  This strategy can indeed be carried out, though it requires a
Vinogradov-type analysis of exponential sums; it also gives the correct asymptotic for \eqref{lar}.  Indeed, this is essentially the approach
taken by van der Corput when establishing infinitely many progressions of length three in the primes.  However, there is a more ``low-tech'' approach
that will give the same qualitative result (but not the asymptotic).  Roughly speaking\footnote{For the detailed rigourous argument, see \cite{green-tao}.}, the idea is as follows.  We allow for the possibility that exponential sums $\E( \Lambda_b(n) e(-\alpha n) | n \in \Z/N\Z )$ could be large, thus providing some additional obstructions to uniformity.
However, the estimate \eqref{lpp} limits the total number of obstructions that could exist.  More precisely, if we introduce a threshold $0 < \eps < 1$
and let $S \subset \Z/N\Z$ denote the exceptional frequencies $\xi$ which obstruct linear uniformity, in the sense that
$$ |\E( \Lambda_b(n) e( - \xi n / N) | n \in \Z/N\Z )| \geq \eps,$$
then \eqref{lpp} shows that $|S| \ll_\eps 1$.  The Vinogradov exponential sum technique will eventually show that $S$ consists only of the zero frequency $0$ for $W,N$ large enough, but we will avoid using this fact, instead treating $S$ as a set for which the only information
known is the cardinality bound.  This approach has the advantage of being more flexible, for instance we will also be able to recover
the result of Green \cite{green} that any subset of the primes with positive relative density contains infinitely many progressions of length three.

The set $S$ represents all the obstructions to uniformity.  We can remove these obstructions by the device of \emph{conditional expectation},
which is a slightly different way than the $W$-trick of removing non-uniformities, though certainly in the same philosophical spirit.  One considers
the \emph{Bohr set} $B(S,\rho) \subset \Z/N\Z$ for some small radius $0 < \rho < 1$ defined by
$$ B(S,\rho) := \{ n \in\Z/N\Z:  \| n\xi \|_{\R/\Z} < \rho \hbox{ for all } \xi \in S \},$$
where $\|x\|_{\R/\Z}$ denotes the distance from $x$ to the nearest integer.
One should think of this Bohr set as being roughly analogous to the subgroup $W\Z$ of $\Z$, thus translates $x+B(S,\rho)$ are the analogues of residue
classes modulo $W$. When executing the $W$-trick, we passed to a single residue class; here, however, we shall proceed in a more ``ergodic'' fashion, averaging out the effect of each translate $x + B(S,\rho)$.  More precisely we split
$$ \Lambda_b = \Lambda_{b,U^\perp} + \Lambda_{b,U}$$ 
where $\Lambda_{b,U^\perp}$ is the ``anti-linearly-uniform'' component
$$ \Lambda_{b,U^\perp}(x) := \Lambda_{b,U^\perp} * \frac{N}{|B(S,\rho)|} 1_{B(S,\rho)} * \frac{N}{|B(S,\rho)|} 1_{B(S,\rho)}(x)$$
where the convolution $f*g$ on $\Z_N$ is defined by
$$ f*g(x) := \E( f(n) g(x-n) | n \in \Z/N\Z),$$
and $\Lambda_{b,U}(x)$ is the ``linearly uniform component''
$$ \Lambda_{b,U} := \Lambda_b - \Lambda_{b,U^\perp}.$$
The function $\Lambda_{b,U^\perp}$ encapsulates all the obstructions to linear uniformity encountered by $\Lambda_b$; the convolution kernel
$$ K := \frac{N}{|B(S,\rho)|} 1_{B(S,\rho)} * \frac{N}{|B(S,\rho)|} 1_{B(S,\rho)}$$
can be thought of as a sort of ``Fej\'er kernel'' adapted to $B(S,\rho)$.  A key observation is that unlike
$\Lambda_b$, the function $\Lambda_{b,U^\perp}$ is \emph{bounded}.  Indeed, from the majorization \eqref{tnub} we have
$$ 0 \leq \Lambda_{b,U^\perp}(x) \ll \tilde \nu_b * K(x)$$
and then by using Fourier expansion of $1_{B(S,\rho)}$ and \eqref{tnub-fourier} one can show
$$ \tilde \nu_b * K(x)
\ll 1 + o_{W \to \infty;|S|,\rho}(1) + o_{N \to \infty;W,|S|,\rho}(1).$$
Since $|S| \ll_\eps 1$, we thus have the uniform boundedness
\begin{equation}\label{lbu}
 0 \leq \Lambda_{b,U^\perp}(x) \ll 1 + o_{W \to \infty;\eps,\rho}(1) + o_{N \to \infty;W,\eps,\rho}(1).
 \end{equation}
In particular we see that $\Lambda_{b,U}$ is pointwise bounded by a constant multiple of $\tilde \nu_b + 1$.
Also, since the kernel $K$ is normalized to have mean $1$, we have
$$ \E( \Lambda_{b,U^\perp}(x) | x \in \Z/N\Z ) = \E( \Lambda_b(x) | x \in \Z/N\Z ) = 1 + o_{W \to \infty}(1) + o_{N \to \infty;W}(1).$$
Thus $\Lambda_{b,U^\perp}$ is bounded, non-negative and has large mean, and so Roth's theorem can be applied (after a renormalization by a bounded scalar) to conclude
\begin{equation}\label{epsrho}
 \E( \Lambda_{b,U^\perp}(a) \Lambda_{b,U^\perp}(a+r) \Lambda_{b,U^\perp}(a+2r) | a,r \in \Z/N\Z ) \geq c - o_{W \to \infty;\eps,\rho}(1) - o_{N \to \infty;W,\eps,\rho}(1)
\end{equation}
for some absolute constant $c > 0$.

The function $\Lambda_{b,U}$ can be regarded as the portion of $\Lambda_b$ remaining after all the obstructions to uniformity have been removed.
By the definition of $S$, one can easily show that $\Lambda_{b,U^\perp}$ has small correlation with all linear phase functions:
$$ |\E( \Lambda_{b,U}(n) e( - \xi n / N) | n \in \Z/N\Z )| \ll \eps + \rho \hbox{ for all } \xi \in \Z/N\Z,$$
and thus by several applications of Lemma \ref{fgh-lin-weighted} we can replace $\Lambda_b$ by $\Lambda_{b,U^\perp}$ with a small error:
\begin{align*}
&\E( \Lambda_b(a) \Lambda_b(a+r) \Lambda_b(a+2r) | a,r \in \Z/N\Z )\\
&\quad\quad= \E( \Lambda_{b,U^\perp}(a) \Lambda_{b,U^\perp}(a+r) \Lambda_{b,U^\perp}(a+2r) | a,r \in \Z/N\Z )
+ O(\eps + \rho).
\end{align*}
Applying \eqref{epsrho} we conclude that
$$ \E( \Lambda_b(a) \Lambda_b(a+r) \Lambda_b(a+2r) | a,r \in \Z/N\Z ) \geq c/2$$
if $\eps,\rho$ are sufficiently small, $W$ is sufficiently large depending on $\eps,\rho$, and $N$ is sufficiently large depending on $\eps,\rho,W$.
This is enough to establish infinitely arithmetic progressions of length threein the primes, and more generally`in any subset of the primes with positive relative
density.  Similar arguments work for other sets that are fairly large and which can be dominated by a suitable enveloping sieve.  For instance,
in \cite{green-tao} it was shown that there were infinitely many arithmetic progressions $p_1,p_2,p_3$ in the primes, where the numbers
$p_1+2$, $p_2+2$, $p_3+2$ are either prime or the product of two primes; this is achieved by combining the arguments above with (a quantitative version of) the famous result of
Chen \cite{Chen} that there are infinitely many primes $p$ such that $p+2$ is the product of at most two primes.

\section{Quadratic obstructions to uniformity}

Let us now consider the task of counting progressions of length four in the primes, or more precisely of obtaining an asymptotic for
$$ \E( \Lambda(a) \Lambda(a+r) \Lambda(a+2r) \Lambda(a+3r) | 1 \leq a,r \leq N ).$$
The Hardy-Littlewood prime tuples conjecture predicts that this quantity is equal to
$\prod_p \alpha_p + o_{N \to \infty}(1)$, where $\alpha_p$ is the local density
$$ \alpha_p := \E( \Lambda_{\Z/p\Z}(a) \Lambda_{\Z/p\Z}(a+r) \Lambda_{\Z/p\Z}(a+2r) \Lambda_{\Z/p\Z}(a+3r) | a,r \in \Z/p\Z ).$$
To put it another way, the number of progressions $a,a+r,a+2r,a+3r$ of primes with $1 \leq a,r \leq N$ is predicted to
be $\frac{N^2}{\log^4 N}(\prod_p \alpha_p + o_{N \to \infty}(1))$.  The result of \cite{gt-primes} establishes a lower bound
$$ \E( \Lambda(a) \Lambda(a+r) \Lambda(a+2r) \Lambda(a+3r) | 1 \leq a,r \leq N ) \geq c - o_{N \to \infty}(1)$$
for some absolute constant $c>0$, which is enough to establish infinitely many progressions of length four in the primes,
but does not give the asymptotic.  In this section we describe a more recent (though significantly more complicated)
approach in \cite{gt-qm}, \cite{gt-mobius}, \cite{gt-inverseu3} which will give the correct asymptotic:

\begin{theorem}\cite{gt-qm}, \cite{gt-mobius}, \cite{gt-inverseu3}  We have
$$ \E( \Lambda(a) \Lambda(a+r) \Lambda(a+2r) \Lambda(a+3r) | 1 \leq a,r \leq N ) = \prod_p \alpha_p + o_{N \to \infty}(1).$$
\end{theorem}

We now sketch the main ideas of proof of this theorem.  Firstly, by the $W$-trick, it will suffice to show that
$$ \E( \Lambda_{b_0}(a) \Lambda_{b_1}(a+r) \Lambda_{b_2}(a+2r) \Lambda_{b_3}(a+3r) | 1 \leq a,r \leq N ) = 1 + o_{W \to \infty}(1) + o_{N \to \infty;W}(1)$$
for all $b_0,\ldots,b_3$ coprime to $W$.  Let us again cheat a little bit by identifying $\{1,\ldots,N\}$ with $\Z/N\Z$ (ignoring some minor truncation issues), so that we now wish to prove that
\begin{equation}\label{lb}
\E( \Lambda_{b_0}(a) \Lambda_{b_1}(a+r) \Lambda_{b_2}(a+2r) \Lambda_{b_3}(a+3r) | a,r \in \Z/N\Z ) = 1 + o_{W \to \infty}(1) + o_{N \to \infty;W}(1).
\end{equation}
It is convenient to take $N$ to be a prime.  We are thus faced with the problem of understanding quartilinear expressions such as
\begin{equation}\label{fghj}
\E( f(a) g(a+r) h(a+2r) j(a+3r) | a,r \in \Z/N\Z );
\end{equation}
to begin the discussion let us suppose that $f,g,h,j$ are bounded in magnitude by $1$.  Let us informally call a function \emph{quadratically uniform} if the above expression is automatically small whenever one of $f,g,h,j$ is replaced with that function.  As in the preceding section, it is easy to see that linear phase functions obstruct quadratic uniformity; however, a new difficulty arises in that \emph{quadratic} phase functions such as
$e(\alpha n^2)$ also obstruct quadratic uniformity.  This can be seen for instance by the identity
\begin{align*}
&\E( f(a) e(-3\alpha(a+r)^2) e(3\alpha(a+2r)^2) e(-\alpha(a+3r)^2) | a,r \in \Z/N\Z ) \\
&\quad = \E( f(n) e(-\alpha n^2) | n \in \Z/N\Z ).
\end{align*}
More generally, one can show that any \emph{quadratic nilsequence} of the form $F(g^n x)$, where $g \in G$ lives in a $2$-step nilpotent Lie group 
$G$, $x$ lives in a compact quotient\footnote{There is an intriguing superficial similarity between the emergence of the $2$-step nilmanifolds $G/\Gamma$ which arise in the analysis of progressions of length $4$, and the cusp manifolds $SL_2(\R) / \Gamma$ which appear for instance in Kloosterman's refinement of the Hardy-Littlewood circle method (which of course corresponds to the unit circle $\R/\Z$).  However, we do not know of a concrete connection between these two different extensions of the circle method.} $G/\Gamma$ of $G$ by a closed subgroup $\Gamma$, and $F: G/\Gamma \to \C$ is a continuous function, will similarly be an obstruction to quadratic uniformity; see \cite{gt-qm}.  The quadratic phases $e(\alpha n^2)$ are good examples of quadratic nilsequence; another example is the generalized quadratic phase $e(\lfloor\alpha n\rfloor \lfloor\beta n\rfloor \gamma)$
for some real numbers $\alpha,\beta,\gamma$, though strictly speakign one needs to smooth out the greatest integer function $\lfloor x \rfloor$ in order
to genuinely obtain a quadratic nilsequence.

The appearance of these quadratic phases shows that the circle method is now insufficient to establish quadratic uniformity; functions such
as $e(\alpha n^2)$ can give significant contributions to \eqref{fghj} despite having very small Fourier coefficients.  However, quadratic uniformity
can still be captured by the very useful \emph{Gowers uniformity norms}\footnote{These are genuine norms for $d \geq 2$; see \cite{gowers}, \cite{gt-primes}, \cite{green-survey}, \cite{tao-survey}.} $U^d(\Z/N\Z)$, defined recursively for $d=0,1,\ldots$ as
$$ \|f\|_{U^0(\Z/N\Z)} := \E( f(x) | x \in \Z/N\Z); \quad \| f \|_{U^{d+1}(\Z/N\Z)} = \E( \| T^h f \overline{f} \|_{U^d(\Z/N\Z)}^{2^d} | h \in \Z/N\Z )^{1/2^{d+1}}$$
where $T^h$ is the shift operator $T^h f(x) := f(x+h)$, thus for instance
\begin{align*}
\| f\|_{U^1(\Z/N\Z)} &= |\E( \overline{f(n)} f(n+h) | n,h \in \Z/N\Z)|^{1/2} \\
&= |\E( f )| \\
\| f\|_{U^2(\Z/N\Z)} &= |\E( f(n) \overline{f(n+h_1) f(n+h_2)} f(n+h_1+h_2) | n,h_1,h_2 \in \Z/N\Z )| \\
&= (\sum_{\xi \in \Z/N\Z} |\hat f(\xi)|^4)^{1/4} \\
\| f\|_{U^3(\Z/N\Z)} &= |\E( \overline{f(n)} f(n+h_1) f(n+h_2) f(n+h_3) \\
&\quad \overline{f(n+h_1+h_2) f(n+h_1+h_3) f(n+h_2+h_3)} f(n+h_1+h_2+h_3) \\
&\quad\quad | n,h_1,h_2,h_3 \in \Z/N\Z )|.
\end{align*}
The relationship between Gowers uniformity norms, and quadratic (or higher order) uniformity, is given by

\begin{lemma}[Generalized von Neumann theorem]\label{gvn}  Let $k \geq 3$, and let $N \geq k-1$ be prime.  If $f_0,\ldots,f_{k-1}: \Z/N\Z \to \C$ are bounded in magnitude by $1$, then
$$ | \E( f_0(a) f_1(a+r) \ldots f_{k-1}(a+(k-1)r) | a,r \in \Z/N\Z ) | \leq \inf_{0 \leq j \leq k} \|f_j\|_{U^{k-1}(\Z/N\Z)}.$$
In particular we have
$$ | \E( f_0(a) f_1(a+r) f_2(a+2r) f_3(a+3r) | a,r \in \Z/N\Z ) | \leq \inf_{0 \leq j \leq 3} \|f_j\|_{U^3(\Z/N\Z)}.$$
\end{lemma}

This lemma can be deduced from $k-1$ applications of the Cauchy-Schwarz inequality, interspersed with $k-1$ applications of the van der Corput identity
$$ |\E( f(n) | n \in \Z/N\Z )|^2 = \E( T^h f(n) \overline{f}(n) | n, h \in \Z/N\Z );$$
we leave the details to the reader (or see \cite{gowers-4}, \cite{gowers}, \cite{host-kra2}, \cite{gt-primes}, \cite{tao:szemeredi}, \cite{green-survey}, \cite{tao-survey}).

The above lemma shows that functions with small $U^3(\Z/N\Z)$ norm are quadratically uniform.  As before, this lemma is not directly applicable to
the problem of finding progressions in primes, since functions such as $\Lambda_b$ are not bounded.  However, because $\Lambda_b$ can be bounded
by an enveloping sieve $\nu_b$ which obeys the good correlation estimates in \eqref{nub}, we can use the following extension of the generalized von Neumann theorem:

\begin{lemma}[Relative generalized von Neumann theorem]\label{gvn-extend}\cite{gt-primes}  Let $k \geq 3$, and let $N > k-1$ be prime.  If $f_0,\ldots,f_{k-1}: \Z/N\Z \to \C$ are such that $f_j$ is bounded by $\nu_{b_j}+1$ for some $b_j$
coprime to $W$, then (if $R = N^{c_k}$ for some sufficiently small $c_k$)
$$ | \E( f_0(a) \ldots f_{k-1}(a+(k-1)r) | a,r \in \Z/N\Z ) | \ll_k \inf_{0 \leq j \leq k} \|f_j\|_{U^{k-1}(\Z/N\Z)} + o_{N \to \infty;W,k}(1)
+ o_{W \to\infty;k}(1).$$
\end{lemma}

This lemma is more complicated to prove than Lemma \ref{gvn} but is still primarily an application of the Cauchy-Schwarz inequality; see\footnote{The argument in \cite{gt-primes} treats the case when all the $b_j$ are equal, but one can easily modify it to treat the case of
distinct $b_j$.}
\cite{gt-primes}, with a heavy reliance on the linear forms estimates \eqref{nub}.  Note that this generalization of Lemma \ref{gvn} is consistent with the transference principle mentioned earlier.

In light of this lemma, we see that in order to establish the asymptotic \eqref{lb}, it will suffice to show that $\Lambda_b - 1$ is quadratically
uniform, or more precisely that
\begin{equation}\label{lb1}
\| \Lambda_b - 1 \|_{U^3(\Z/N\Z)} = o_{N \to \infty;W}(1)
+ o_{W \to\infty}(1)
\end{equation}
for all $b$ coprime to $W$.  This is not easy to do directly, since the quantity $\| \Lambda_b - 1 \|_{U^3(\Z/N\Z)}$ is basically the same
type of expression that appears in the Hardy-Littlewood prime tuples conjecture, and is beyond the reach of the circle method.  Nevertheless, one
can proceed by locating all the obstructions to quadratic uniformity, and then checking that the function $\Lambda_b - 1$ is orthogonal to all of these.

We have already observed that the quadratic nilsequences $F(g^n x)$ are obstructions to quadratic uniformity.  Recent developments \cite{host-kra2}, \cite{bhk} in ergodic theory strongly suggest\footnote{Roughly speaking, the ergodic theory setting corresponds to considering
averages such as $\E( f(a) f(a+r) f(a+2r) f(a+3r) | 1 \leq a \leq N, 1 \leq r \leq H )$ where the shift range $H$ goes to infinity much more slowly than $N$ does.  As such, there does not appear to be a direct ``correspondence principle'' between the results in \cite{host-kra2}, \cite{bhk} and the type of results considered here, but there is certainly a very strong analogy between the two.  See \cite{kra-survey} for more on the ergodic theory perspective to these problems.} that these are in fact the only obstructions to quadratic uniformity.  By building on the pioneering combinatorial and analytical technology of Gowers \cite{gowers-4}, a quantitative version of this assertion was made in \cite{gt-qm}.  More precisely:

\begin{theorem}[Inverse theorem for $U^3(\Z/N\Z)$]\label{gt-inverse}\cite{gt-qm}  Let $0 < \eta < 1$.  Then there exists a collection ${\mathcal N}$ of $O_\eta(1)$ triples $(G,\Gamma,F)$, where
$G$ is a $2$-step nilpotent Lie group, $\Gamma$ is a closed co-compact subgroup of $G$, and $F: G/\Gamma \to \C$ is a smooth function, with the
following property: if $N$ is an odd prime and $f: \Z/N\Z \to \C$ is bounded by $1$ and is such that $\|f\|_{U^3(\Z/N\Z)}$, then there exists
a triple $(G,\Gamma,F)$ from this collection, a group element $g \in G$, a point $x \in G/\Gamma$, and a shift $h \in \Z/N\Z$ such that
$$ |\E( T^h f(n) \overline{F(g^n x)} | -N/2 < n < N/2 )| \gg_\eta 1.$$
\end{theorem}

One can explicitly describe the collection ${\mathcal N}$, and give quantitative bounds on the dimension of $G/\Gamma$ and the smoothness of $F$, as
well as the dependence of the implied constant on $\eta$; see \cite{gt-qm}.

The proof of Theorem \ref{gt-inverse} is quite lengthy, using many tools of Gowers in additive combinatorics and Fourier analysis.  On the
other hand, it may well be that a ``softer'' proof, without the quantitative bounds, is available by the ergodic-theory methods in \cite{host-kra2}, \cite{bhk}.  In \cite{gt-inverseu3}, the results from \cite{gt-primes} (and more precisely, Theorem \ref{structure} below) 
were used to extend Theorem \ref{gt-inverse} to the case when $f$ is merely bounded by $\nu_b + 1$ rather than by $1$; again, this is consistent
with the transference principle.  By applying this extended version of Theorem \ref{gt-inverse}, we see that one can prove \eqref{lb1} as soon as
one demonstrates the asymptotic orthogonality estimate
\begin{equation}\label{lambda-ortho}
 \E( (T^h \Lambda_b(n) - 1) \overline{F(g^n x)} | -N/2 < n < N/2 ) = o_{N \to \infty;W,F,G,\Gamma}(1) + o_{W \to\infty;F,G,\Gamma}(1)
 \end{equation}
for all quadratic nilsequences $F(g^n x)$.  

This type of result is essentially an exponential sum estimate on $\Lambda$, and can thus be attacked by the standard Vinogradov-type methods.  A model case is the estimate
$$ \E( (\Lambda_b(n)-1) e(-\alpha n^2) | 1 \leq n \leq N ) = o_{N \to \infty}(1)$$
for all $\alpha \in \R$, which was essentially obtained in \cite{ghosh}.  The general case of quadratic nilsequences is treated in \cite{gt-mobius},
\cite{gt-inverseu3}.  In those papers it is convenient to first prove the preliminary estimate
$$ \E( \mu(n) \overline{F(g^n x)} | 1 \leq n\leq N ) \ll_{A,F,G,\Gamma} \log^{-A} N$$
for all $A > 0$ whenever $F$ is smooth; see \cite{gt-inverseu3}.  This can be considered a generalization of 
Davenport's estimate\cite{davenport}
$$ \E( \mu(n) e(-\alpha n) | 1\leq n \leq N ) \ll_A \log^{-A} N$$
and is proven by broadly similar, though significantly more technical, methods (in particular, Vaughan's identity, a division into major and minor arcs, and Cauchy-Schwarz type arguments to deal with the minor arcs).  It is however simpler to deal with the M\"obius function $\mu(n)$ than
the modified von Mangoldt function $\Lambda_b(n) - 1$, as $\mu$ is bounded, and also obeys a somewhat more pleasant Vaughan identity than
$\Lambda$.  Using this estimate and some elementary arguments, 
it is already possible to establish
$$ \E( (T^h \Lambda_b(n) - T^h \Lambda_{R,\varphi,b}(n)) \overline{F(g^n x)} | -N/2 < n < N/2 ) = o_{N \to \infty;W,F,G,\Gamma}(1) + o_{W \to\infty;F,G,\Gamma}(1)$$
where $\Lambda_{R,\varphi,b}(n) := \frac{W}{\phi(W)} \Lambda_{R,\varphi}(Wn+b)$ and $\Lambda_{R,\varphi}$ was defined\footnote{Actually, any reasonable truncated divisor sum approximation to $\Lambda$ could be used in place of $\Lambda_{R,\varphi}$ here.} 
in \eqref{lurch}; as usual we set $R$ to be a small power of $N$ and $\varphi$ to be a suitable cutoff function.  By the triangle inequality, it thus remains to verify that
$$ \E( (T^h \Lambda_{R,\varphi,b}(n) - 1) \overline{F(g^n x)} | -N/2 < n < N/2 ) = o_{N \to \infty;W,F,G,\Gamma}(1) + o_{W \to\infty;F,G,\Gamma}(1).$$
It turns out that the simplest way to do this is to apply the Cauchy-Schwarz inequality (in the spirit of Lemma \ref{gvn} and Lemma \ref{gvn-extend},
and in particular on the \emph{Gowers-Cauchy-Schwarz inequality} introduced in \cite{gowers}, and also playing a key role in \cite{gt-primes}),
to reduce matters to the $U^3$ estimate
$$ \| \Lambda_{R,\varphi,b}(n) - 1 \|_{U^3(\Z/N\Z)} = o_{N \to \infty;W}(1) + o_{W \to\infty}(1),$$
which in turn can be established by a Goldston-Y{\i}ld{\i}r{\i}m correlation estimate, similar in spirit to \eqref{nub}.  See \cite{gt-inverseu3}.

It is entirely possible that the techniques discussed in this section extend to give an asymptotic for longer progressions in the primes, 
though there are serious new difficulties that appear (similar to the new difficulties that appear in \cite{gowers} when compared against \cite{gowers-4}).  We (in joint work with Ben Green) hope to report on this problem in a future paper.

\section{Ergodic obstructions to uniformity}

In the previous section, we outlined a rather complicated approach that yielded an asymptotic for the number of progressions of length four
in the primes.  As we already saw though in the length three case, it can often be significantly easier to establish the weaker result of a
non-trivial lower bound for the number of such progressions, using tools such as Roth's theorem.  This was achieved in \cite{gt-primes}, in particular
establishing that the primes contain arbitrarily long arithmetic progressions.  The argument can be seen as a variant of the above arguments,
but in which the ``hard'' obstructions of nilsequences are replaced by much ``softer'' obstructions coming from ergodic averages.  These soft obstructions are insufficiently explicit to easily allow for establishing asymptotic orthogonality results such as \eqref{lambda-ortho}, but
they are still controllable to the extent that one can modify the arguments of Section \ref{fourier-sec}, using the soft obstructions to build
generalized Bohr sets with which to split $\Lambda_b$ into a uniform component, which is negligible, and an anti-uniform component, which can
be treated by a theorem of Szemer\'edi.

We turn to the details. The famous theorem of Szemer\'edi \cite{szemeredi} asserts that every subset of integers of positive density contains
arbitrarily long arithmetic progressions.  A quantitative version of this theorem, which generalizes Theorem \ref{qrt}, is as follows:

\begin{theorem}[Quantitative Szemer\'edi theorem]\label{qst} Let $k \geq 1$, and let $f: \Z/N\Z \to \R$ be a function such that $0 \leq f(n) \leq 1$ for all $n \in \Z/N\Z$, and such that $\E( f(n) | n \in \Z/N\Z ) \geq \delta$ for some $0 < \delta < 1$.  Then we have
$$ \E( f(a) f(a+r) \ldots f(a+(k-1)r) | a,r \in \Z/N\Z ) \geq c(k,\delta)$$
for some $c(k,\delta) > 0$.
\end{theorem}

This theorem can be deduced from Szemer\'edi's original theorem from the averaging argument of Varnavides \cite{varnavides}; see also
\cite{tao:szemeredi} for a direct proof.

As in Section \ref{fourier-sec}, the task (after applying the $W$-trick) is to obtain a non-trivial lower bound for
\begin{equation}\label{labb}
\E( \Lambda_b(a) \ldots \Lambda_b(a+(k-1)r) | a,r\in \Z/N\Z),
\end{equation}
where we once again gloss over the distinction between $\Z/N\Z$ and $\{1,\ldots,N\}$ to simplify the discussion.  Again, we cannot apply
Theorem \ref{qst} directly because of the unboundedness of $\Lambda_b$.  However, we can proceed by establishing the following structure theorem, that decomposes any non-negative function bounded by the enveloping sieve $\nu_b$ into a Gowers uniform component (with small Gowers uniformity norm), a non-negative bounded component, and a small error.

\begin{theorem}[Structure theorem]\label{structure}\cite{gt-primes}  Let $k \geq 1$, and let $R = N^{c_k}$ for some sufficiently 
small $c_k > 0$.  Let $f: \Z/N\Z \to \R$
be such that $0 \leq f(n) \leq \nu_b(n)$.  Let $0 < \eps < 1$.  Then functions $f_U, f_{U^\perp}: \Z/N\Z \to \C$
such that
\begin{equation}\label{fuu}
 \| f_U \|_{U^{k-1}(\Z/N\Z)} = o_{\eps \to 0;k}(1)
\end{equation}
and
\begin{equation}\label{fup}
 0 \leq f_{U^\perp}(n) \leq 1 + o_{\eps \to 0;k}(1) + o_{N \to \infty;\eps,k}(1)
 \end{equation}
and
$$ 0 \leq f_U(n) + f_{U^\perp}(n) \leq f(n)$$
for all $n \in \Z/N\Z$.  Furthermore, we have 
\begin{equation}\label{fun-mean-diff}
\E( |f(n) - f_{U^\perp}(n) - f_U(n)| | n \in \Z/N\Z ) = o_{\eps \to 0;k}(1).
\end{equation}
and
\begin{equation}\label{fun-mean}
\E( f_{U^\perp}(n) | n \in \Z/N\Z ) = \E( f(n) | n \in \Z/N\Z ) + o_{\eps \to 0;k}(1).
\end{equation}
\end{theorem}

Assuming this theorem, a lower bound for \eqref{labb} can be easily accomplished.  By \eqref{tnub} we can apply Theorem \ref{structure} with $f := c \Lambda_b$ for some absolute constant $c > 0$, to obtain a majorization
$$ 0 \leq f_U + f_{U^\perp} \leq c \Lambda_b.$$
It then suffices to obtain a lower bound for
$$
\E( (f_U + f_{U^\perp})(a) \ldots (f_U + f_{U^\perp})(a+(k-1)r) | a,r\in \Z/N\Z).
$$
All the terms involving at least one factor of $f_U$ are $o_{\eps \to 0;k}(1) + o_{N \to \infty;\eps,k}(1)$, thanks mainly to \eqref{fuu} and Lemma \ref{gvn-extend}.  The remaining term involving $f_{U^\perp}$ is at least $c_k - o_{\eps \to 0; k}(1) - o_{N \to \infty;\eps,k}(1)$,
thanks to Theorem \ref{qst} and \eqref{fun-mean}.  Setting $\eps$ suitably small, and then $N$ sufficiently large, we obtain a non-trivial lower bound for \eqref{labb}.

Thus Theorem \ref{structure} allows one to transfer Theorem \ref{qst} to a relative setting, adapted to the enveloping sieve $\nu_b$.  A similar
argument also allows one to use Theorem \ref{structure} to transfer Theorem \ref{gvn-extend} to the relative setting; see \cite{gt-inverseu3}.

It remains to prove Theorem \ref{structure}.  Let us fix $f$.  The first guess is to take $f_{U^\perp}$ to be the mean of $f$,
$f_{U^\perp} := \E(f)$, and then set $f_U := f - f_{U^\perp}$.  It is clear that $f_{U^\perp}$ is non-negative, and also
$$ f_{U^\perp} = \E(f) \leq \E(\nu_b) = 1 + o_{W \to \infty;k}(1) + o_{N \to \infty;W,k}(1).$$
Also we trivially have \eqref{fun-mean} and \eqref{fun-mean-diff}.  The only difficulty is that we do not necessarily have \eqref{fuu}; there is
no reason why $f_U$ needs to be Gowers uniform (i.e. have small $U^{k-1}(\Z/N\Z)$ norm).  However, if this is the case, it turns out to be
possible to locate a precise obstruction which is preventing $f_U$ from being uniform, and transfer this obstruction from $f_{U}$ to $f_{U^\perp}$.
This may not remove all the non-uniformity from $f_U$, but it will increase the energy ($L^2(\Z/N\Z)$ norm) of $f_{U^\perp}$ by a significant
amount, and so after iterating this process a finite number of times we will eventually end up with a Gowers uniform $f_U$.  

The above type of argument has also been used before in ergodic theory (most notably in Furstenberg's structure theorem \cite{FKO}), and also
in the proof of the Szemer\'edi regularity lemma \cite{szemeredi}; not co-incidentally, both of those cited papers concerned Szemer\'edi's theorem
(Theorem \ref{qst}).  The argument in Section \ref{fourier-sec} involving convolution with a Bohr set generated by all the Fourier obstructions
to uniformity is also an argument of this type (although in that case one transferred all the obstructions from $f_U$ to $f_{U^\perp}$ at once, rather than one at a time).  The main difficulty in executing the above idea is to maintain \eqref{fup} throughout this procedure, i.e. to 
keep $f_{U^\perp}$ non-negative and bounded by $1$ (plus negligible errors).  To achieve the non-negativity, the simplest way is to use
the machinery of conditional expectation (as is done in Furstenberg's structure theorem, and implicitly in the Szemer\'edi regularity lemma).
To achieve the boundedness, one needs some control on the obstructions to uniformity that one is transferring to $f_{U^\perp}$.  In the Fourier-analytic argument, these obstructions are linear phase functions $e(\alpha n)$, and one can use Fourier-analytic control in the enveloping
sieve (see \eqref{tnub-fourier}) to keep $f_{U^\perp}$ bounded.  To adopt a similar argument in the general case, one might imagine one would need a similarly explicit description of these obstructions, for instance using the nilsequences of the preceding section.  However, it turns out that one can get by using a much less explicit obstruction to uniformity, first introduced in ergodic theory\footnote{More precisely, the key observation for ergodic theory is that the obstructions to weak mixing (which roughly corresponds to Gowers uniformity) are given by almost periodic functions, and more specifically given any function $f$ which fails to be weakly mixing (so that $\langle T^h f, f \rangle$ does not converge on average to zero), one can construct the non-trivial almost periodic function $F := \lim_{H \to \infty} \E( \langle T^h f, f \rangle T^h f | -H \leq h \leq H )$, which has a positive correlation with $f$.  See for instance \cite{FKO}; for the connection with the Gowers uniformity norms see \cite{host-kra2}, \cite{kra-survey}.}.

In order to make the above strategy rigourous, we need two basic concepts, that of a \emph{dual function} and that of \emph{conditional expectation}.
The dual function ${\mathcal D}_d f: \Z/N\Z \to \C$ of a function $f: \Z/N\Z \to \C$ is defined recursively for $d = 0,1,2,\ldots$ by the formula
$$ {\mathcal D}_0 f = 1; \quad {\mathcal D}_{d+1} f = \E( {\mathcal D}_d(f \overline{T^h f}) T^h f | h \in \Z/N\Z );$$
thus for instance 
\begin{align*}
{\mathcal D}_1 f(n) &= \E(f) \\
{\mathcal D}_2 f(n) &= \E( f(n+h_1) f(n+h_2) \overline{f(n+h_1+h_2)} | h_1,h_2 \in \Z/N\Z) \\
&= \E( \langle f, T^h f \rangle T^h f(n) | h \in \Z/N\Z ) \\
&= \sum_{\xi \in \Z/N\Z} |\hat f(\xi)|^2 f(\xi) e(n \xi/N) \\
{\mathcal D}_3 f(n) &= \E( f(n+h_1) f(n+h_2) f(n+h_3) \overline{f(n+h_1+h_2) f(n+h_1+h_3) f(n+h_2+h_3)} \\
&\quad f(n+h_1+h_2+h_3) | h_1,h_2,h_3 \in \Z/N\Z)
\end{align*}
where $\langle,\rangle$ denotes the usual inner product $\langle f, g \rangle = \E(f \overline{g})$.
One can easily use induction to verify that
\begin{equation}\label{fkf}
\langle f, {\mathcal D}_{k-1} f \rangle = \|f\|_{U^{k-1}(\Z/N\Z)}^{2^{k-1}}.
\end{equation}
Thus if $f$ fails to be Gowers uniform of order $k-1$, it correlates with a dual function ${\mathcal D}_{k-1} f$.  These dual functions will serve
as our obstructions to Gowers uniformity; they are simple to describe but are not very explicit, as they involve a function $f$ for which we have only limited control.  Nevertheless, there is a large amount of averaging contained in the non-linear operator ${\mathcal D}_{k-1}$, which will allow us to obtain satisfactory control on these dual functions.

To proceed further, we need to understand the properties of dual functions better.  The first important (and easy) property is that dual functions are always bounded:
more precisely, we have $|{\mathcal D}_{k-1} f| \ll_k 1$ whenever $f$ is pointwise bounded by $\nu_b + 1$.  Indeed, in such a case we have
$$|{\mathcal D}_{k-1} f| \leq {\mathcal D}_{k-1}(\nu_b+1),$$
and several applications of \eqref{nub} gives the bound ${\mathcal D}_{k-1}(\nu_b+1)$ (see \cite{gt-primes}).

The second important (but significantly deeper) property is that a dual function, and more generally any polynomial combination of dual functions, is highly ``Gowers anti-uniform'' in the sense that it is essentially orthogonal to all Gowers uniform functions, and in particular to the function $\nu_b - 1$ (which can easily be shown to be Gowers uniform, thanks to several applications of \eqref{nub}).  Indeed, it turns out that we have
\begin{equation}\label{scat}
 \langle \nu_b - 1, P( {\mathcal D}_{k-1}(f_1), \ldots {\mathcal D}_{k-1}(f_m) ) \rangle = o_{N \to \infty; m,P,W}(1) + o_{W \to \infty; m,P}(1)
\end{equation}
for any polynomial $P(x_1,\ldots,x_m)$ of $m$ variables, and any functions $f_1,\ldots,f_m: \Z/N\Z \to \C$ bounded in magnitude by $\nu_b + 1$.
This fact is elementary to prove, but not entirely trivial; it is obtained by a large number of applications of the Cauchy-Schwarz and H\"older inequalities, combined with the correlation condition \eqref{corrtau}.  See \cite{gt-primes}.

One should compare the above facts with the situation in the Fourier-analytic argument.  In that argument, the role of dual functions was played by the linear phase functions $e(\alpha n)$, which are certainly bounded.  A polynomial combination of linear phase functions is nothing more than a trigonometric polynomial, and \eqref{tnub-fourier} then shows that $\nu-1$ is indeed mostly orthogonal to such polynomial combinations.

To exploit these facts about dual functions, we need to introduce the machinery of $\sigma$-algebras and conditional expectation.

\begin{definition}  A \emph{$\sigma$-algebra} is a collection ${\mathcal B}$ of subsets of $\Z/N\Z$ which contains $\emptyset$ and $\Z/N\Z$ and is closed under union, intersection, and complementation.  A function $f: \Z/N\Z \to \C$ is \emph{${\mathcal B}$-measurable} if all its level sets lie in ${\mathcal B}$.
If ${\mathcal B}$ is a $\sigma$-algebra and $f: \Z/N\Z \to \C$, we define the \emph{conditional expectation}
$\E(f|\B): \Z/N\Z \to \C$ of $f$ with respect to $\B$ to be the function
$$ \E(f|\B)(x) := \E(f|\B(x)) = \frac{1}{\B(x)} \sum_{n \in \B(x)} f(n)$$
for all $x \in \Z/N\Z$ where $\B(x)$ is the smallest set in $\B$ which contains $x$.
If $\B_1, \B_2$ are two $\sigma$-algebras, we use $\B_1 \vee \B_2$ to denote the smallest $\sigma$-algebra which contains both $\B_1$ and $\B_2$.
\end{definition}

A basic fact in measure theory is that any algebra of functions generates a $\sigma$-algebra.  The estimate \eqref{scat} asserts, morally speaking,
 that $\nu_b-1$ is asymptotically orthogonal to the algebra generated by dual functions, and thus should also be orthogonal to the $\sigma$-algebra generated by dual functions.  Indeed, we can make this precise as follows.  Given any dual function ${\mathcal D}_{k-1}(f)$ and any cutoff $\eps > 0$, we can generate a $\sigma$-algebra $\B_\eps({\mathcal D}_{k-1}(f))$, by partitioning the complex plane $\C$ into squares of length $\eps$, and using the inverse images of these squares under ${\mathcal D}_{k-1}(f)$ as the atoms of the $\sigma$-algebra.  There is some choice in how to choose this partition; a random translation of the standard partition will work here.  A key result in \cite{gt-primes} is then that for any $m \geq 1$ and any functions
 $f_1,\ldots,f_m$ bounded in magnitude by $\nu_b + 1$, we have the uniform distribution property
\begin{equation}\label{nubile}
\E( \nu_b-1 | \B_\eps({\mathcal D}_{k-1}(f_1)) \vee \ldots \vee \B_\eps({\mathcal D}_{k-1}(f_m)) = o_{\eps \to 0;m,k}(1) + o_{W \to \infty;m,k,\eps}(1)
+ o_{N \to \infty;m,k,\eps,N}(1)
\end{equation}
except on an exceptional set $\Omega$ which is small in the sense that
$$ \E( (\nu_b + 1) 1_{\Omega} ) = o_{\eps \to 0;m,k}(1) + o_{W \to \infty;m,k,\eps}(1)
+ o_{N \to \infty;m,k,\eps,N}(1).$$
This claim can be derived fairly quickly from \eqref{scat} and the Weierstrass approximation theorem\footnote{As our functions here are complex valued, we have to consider polynomials which involve the conjugates of the dual functions ${\mathcal D}_{k-1}(f_j)$ as well as the dual functions themselves,
but this does not cause any additional difficulty}; see \cite{gt-primes}.

We can now sketch the proof of Theorem \ref{structure}.  As mentioned earlier, the idea is to detect any obstructions to uniformity in $f_U$
(in the guise of dual functions ${\mathcal D}_{k-1}(f_1), \ldots, {\mathcal D}_{k-1}(f_m)$, where $f_1,\ldots,f_m$ are bounded in magnitude by $\nu_b + 1$) and transfer them to $f_{U^\perp}$ one at a time.  Oversimplifying somewhat (in particular, glossing over the role of the exceptional set $\Omega$), the algorithm for doing so is as follows:

\begin{itemize}

\item Step 0.  Set $m=0$.

\item Step 1.  Set $f_{U^\perp} := \E( f | \B_\eps({\mathcal D}_{k-1}(f_1)) \vee \ldots \vee \B_\eps({\mathcal D}_{k-1}(f_m))$ (so initially we would
just have $f_{U^\perp} = \E(f)$), and then set $f_U := f - f_{U^\perp}$.  Clearly $f_{U^\perp}$ is non-negative and has the same mean as $f$; 
from \eqref{nubile} we ensure that $f_{U^\perp}$ is bounded.

\item Step 2.  If $f_U$ is Gowers uniform, in the sense that $\|f_U \|_{U^{k-1}(\Z/N\Z)} \leq \eps^{1/2}$, then we are done.  Otherwise, we set $f_{m+1} := f_U$, increment $m$ by $1$, and return to Step 1.

\end{itemize}

It turns out that every time we return from Step 2 to Step 1, the energy $\E( |f_{U^\perp}|^2 )$ of $f_{U^\perp}$ increases by at least
$c_{\eps,k}$ (plus some negligible $o(1)$ errors), where $c_{\eps,k} > 0$ is an explicit positive quantity depending only on $\eps$ and $k$; see
\cite{gt-primes}.  Intuitively, the reason for this is as follows.  If $f_U$ is not Gowers uniform, then by \eqref{fkf} $f_U$ as a large correlation with ${\mathcal D}_{k-1}(f_U) = {\mathcal D}_{k-1}(f_{m+1})$.  But $f_U$, by construction, is orthogonal to all the functions which are measurable with respect to the $\sigma$-algebra $\B_\eps({\mathcal D}_{k-1}(f_1)) \vee \ldots \vee \B_\eps({\mathcal D}_{k-1}(f_m)$, while ${\mathcal D}_{k-1}(f_{m+1})$
lies (modulo negligible errors) in the larger $\sigma$-algebra $\B_\eps({\mathcal D}_{k-1}(f_1)) \vee \ldots \vee \B_\eps({\mathcal D}_{k-1}(f_{m+1})$.
The energy increment then follows (morally, at least) from the following simple lemma:

\begin{lemma}[Correlation implies energy increment]  Let $\B \subseteq \B'$ be $\sigma$-algebras, and let $f, g$ be functions such that $f$ is orthogonal to all $\B$-measurable functions, while $g$ is $\B'$-measurable and bounded in magnitude by $1$.  Then we have the energy increment
$$ \E( |\E(f|\B')|^2 ) \geq \E( |\E(f|\B)|^2 ) + |\langle f, g \rangle|^2.$$
\end{lemma}

\begin{proof}  From the $\B'$-measurability of $g$ we have
$$ \langle f, g \rangle = \langle \E(f|\B'), g \rangle.$$
Also, since $f$ is orthogonal to all $\B$-measurable functions, we have $\E(f|\B) = 0$.  Thus
$$ \langle f, g \rangle = \langle \E(f|\B') - \E(f|\B), g \rangle.$$
Applying Cauchy-Schwarz and the boundedness of $g$ we conclude
$$ \E( |\E(f|\B') - \E(f|\B)|^2 ) \geq |\langle f, g \rangle|^2$$
and the claim then follows from Pythagoras' theorem.
\end{proof}

In practice, we cannot quite use this simple lemma because of the presence of the exceptional sets
$\Omega$, but it is still possible to obtain the energy increment by carefully modifying the above argument; see \cite{gt-primes}.

Observe that the energy  $\E( |f_{U^\perp}|^2 )$ increments by a fixed factor at each stage of the iteration, but remains bounded independently of
the number of steps of the iteration (ignoring some negligible $o(1)$ type errors).  Thus the algorithm can only run for a bounded number of steps,
which keeps all the $o(1)$ errors under control.  After doing all the book-keeping, one eventually arrives at a proof of Theorem \ref{structure};
see \cite{gt-primes} for the full details.  As discussed earlier, this is enough to establish that the primes contain arbitrarily long arithmetic progressions; the same argument also shows that any subset of the primes of positive relative density contain arbitrarily long arithmetic progressions.
One can also follow through the argument carefully to eventually yield a lower bound
$$ \E( \Lambda(a) \ldots \Lambda(a+(k-1)r) | 1 \leq a,r \leq N ) \geq c(k) - o_{N \to \infty;k}(1)$$
for some explicitly computable $c(k) > 0$; the exact value is rather poor, depending on both the quantitative error bounds in 
the correlation estimates \eqref{nub}, \eqref{corrtau}, as well the constant in Theorem \ref{qst}.

\section{Further directions}

The transference methods here should be applicable to some other situations.  For instance, a variant of the above argument
was used recently in \cite{tao:gaussian} to show that the Gaussian primes in $\Z[i]$ contain infinitely many constellations of any prescribed
shape and orientation; one needs to replace Szemer\'edi's theorem by the somewhat stronger ``hypergraph removal lemma'' of Gowers \cite{gowers-reg}
and R\"odl-Skokan \cite{rodl}, \cite{rodl2} (see also \cite{tao:hyper}), and the presence of the conjugation
operation $z \mapsto \overline{z}$ in the Galois group $Gal(\Q[i]/\Q)$ causes some technical difficulties,
but otherwise the strategy is almost identical.  We refer the reader to \cite{tao:hyper} and \cite{tao:gaussian} for further details.
Similar results should also hold for other number fields that enjoy unique factorization.  For instance, one should be able to show that given any
finite field $F$, the monic irreducible polynomials of one variable in $F[x]$ should contain affine subspaces over $F$ of arbitrarily high dimension.

A more challenging extension would be to obtain a multidimensional relative Szemer\'edi theorem, which would assert that
given any dimension $d \geq 1$, and given the set of primes $P = \{2,3,5,\ldots\}$, that any subset of $P^d$ of positive relative density should contain
infinitely many constellations of any prescribed shape and orientation.  For $P^d$ replaced by $\R^d$, this result was proven in
\cite{fk}, and also follows from the hypergraph removal lemma mentioned briefly earlier.  A major new difficulty here is that the natural enveloping
sieve for $P^d$ is not very pseudorandom, even after applying the higher-dimensional analogue of the $W$-trick; the lack of pseudorandomness, even for $P^2$, can be seen by the observation that if the two acute corners of a right-angled triangle (with sides parallel to the axes) lie in $P^2$, then
the third corner also automatically lies in $P^2$, despite $P^2$ being quite sparse.  We do not know how to resolve this problem.

It should also be possible to establish arbitrarily long progressions $a, a+r,\ldots, a+(k-1)r$ in the primes (or any positive relative density subset
thereof), in which the spacing $r$ is significantly smaller than the base point $a$, obtaining for instance progressions such that $r = O_{\eps,k}(a^\eps)$ for any given $\eps$.  This is likely to follow by localizing the above theory to intervals of length $O(N^{\eps})$
in $\{N+1,\ldots,2N\}$.

A more difficult result would be to obtain a polynomial Szemer\'edi theorem for the primes. More precisely, if $P_1,\ldots,P_k: \Z \to \Z$
were any polynomials mapping the integers to the integers with $P_1(0) = \ldots = P_k(0) = 0$, then there should be infinitely many $k$-tuplets 
$a + P_1(r),\ldots,a+P_k(r)$ with $r \neq 0$, such that all the $a + P_j(r)$ are prime.  If the primes were replaced by a positive density subset
of $\Z$, then this result was obtained by Bergelson and Leibman \cite{bergelson-leibman}.  If one wished to localize this problem to $\Z/N\Z$, it would
be necessary to restrict $r$ to be at most a small power of $N$, and so one may first have to understand the previous  problem concerning
progressions with small spacing before tackling this problem. The hypothesis $P_1(0) = \ldots = P_k(0)= 0$ seems to unfortunately be rather crucial
to the method (for instance, one can easily construct counterexamples to the Bergelson-Leibman theorem without this hypothesis), which is a pity
as one would otherwise have a route to prove such conjectures as the twin primes conjecture or more generally the Hardy-Littlewood prime tuple
conjecture.

Another problem (communicated by Vitaly Bergelson) which might now be feasible is to establish that the set $P-1 = \{ 1, 2, 4, 6, 10, \ldots\}$ formed by decrementing one from each prime, is an IP set, or more precisely given any $k$ there exist distinct $a_1,\ldots,a_k$ such that the finite sums
$\{ \sum_{j \in J} a_j: J \subseteq \{1,\ldots,k\}; J \neq \emptyset \}$ are contained in $P-1$.  The case $k=2$ can be handled by the circle method,
but the higher $k$ remain open.  Such a result would then lead to a number of combinatorial consequences, see for instance \cite{br} for further discussion.

\providecommand{\MR}{\relax\ifhmode\unskip\space\fi MR }
\providecommand{\MRhref}[2]{%
  \href{http://www.ams.org/mathscinet-getitem?mr=#1}{#2}
}
\providecommand{\href}[2]{#2}

     \end{document}